\documentclass[a4paper]{scrartcl}
       
\usepackage[utf8]{inputenc}
\usepackage[T1]{fontenc}
\usepackage[english]{babel} 
\usepackage{lmodern}
\usepackage[]{algorithm2e}
\usepackage[onehalfspacing]{setspace}
\usepackage{amsmath,amssymb,amsthm,amsfonts,amsbsy,latexsym}
\usepackage{graphicx}
\usepackage{geometry}
\usepackage{natbib}
\usepackage{multirow}
\usepackage{float}
\usepackage{subfloat}
\usepackage{subfig}
\usepackage{fancyhdr}
\usepackage{makecell}
\usepackage{pgf,tikz}
\usepackage{makeidx}
\usetikzlibrary{shapes,snakes}
\usetikzlibrary{fadings}
\usepackage[pdfpagelabels=true]{hyperref}
\usepackage{tocbibind}

\newcommand{\R}{\mathbb{R}}

\newcommand{\N}{\mathbb{N}}

\DeclareMathOperator*{\argmin}{argmin}
\DeclareMathOperator*{\argmax}{argmax}

\newtheorem{thm}{Theorem}[section]
\newtheorem{Def}{Definition}[section]

\newtheorem{rem}{Remark}[section]
\newtheorem{prop}{Proposition}[section]
\newtheorem{cor}{Corollary}[section]
\newtheorem{ex}{Example}[section]

\newtheorem*{bew}{Proof}

\numberwithin{equation}{section}

\begin{document}

\title{Trimming Stability Selection increases variable selection robustness}
\author{Tino Werner\footnote{Institute for Mathematics, Carl von Ossietzky University Oldenburg, P/O Box 5634, 26046 Oldenburg (Oldb), Germany, \texttt{tino.werner1@uni-oldenburg.de}}}
\maketitle

\begin{footnotesize} 

\begin{abstract} Contamination can severely distort an estimator unless the estimation procedure is suitably robust. This is a well-known issue and has been addressed in Robust Statistics, however, the relation of contamination and distorted variable selection has been rarely considered in literature. As for variable selection, many methods for sparse model selection have been proposed, including the Stability Selection which is a meta-algorithm based on some variable selection algorithm in order to immunize against particular data configurations. We introduce the variable selection breakdown point that quantifies the number of cases resp. cells that have to be contaminated in order to let no relevant variable be detected. We show that particular outlier configurations can completely mislead model selection and argue why even cell-wise robust methods cannot fix this problem. We combine the variable selection breakdown point with resampling, resulting in the Stability Selection breakdown point that quantifies the robustness of Stability Selection. We propose a trimmed Stability Selection which only aggregates the models with the lowest in-sample losses so that, heuristically, models computed on heavily contaminated resamples should be trimmed away. An extensive simulation study with non-robust regression and classification algorithms as well as with Sparse Least Trimmed Squares reveals both the potential of our approach to boost the model selection robustness as well as the fragility of variable selection using non-robust algorithms, even for an extremely small cell-wise contamination rate. \end{abstract} \end{footnotesize}

\section{Introduction}

Feature selection in machine learning became popular once Tibshirani introduced the Lasso (\cite{tibsh96}). This opened the path for a plethora of feature selection methods in regression (e.g., \cite{bu03}, \cite{efron04}, \cite{zou06}, \cite{yuanlin}, \cite{hastie13}), classification (\cite{park07}, \cite{mein07}, \cite{geer08}), clustering (e.g., \cite{tibsh10}, \cite{alelyani} and references therein), ranking (\cite{tian}, \cite{lai13}, \cite{laporte}) and sparse covariance or precision matrix estimation (\cite{banerjee08}, \cite{friedman08}, \cite{geer16}). 

It has been observed that single regularized models often tend to overfit which was the starting point for a more sophisticated concept that combines ensemble models with feature selection, namely Stability Selection (\cite{bu10}). Roughly speaking, in Stability Selection one considers subsamples of the training data, performs feature selection on each subsample and aggregates the models to get a final stable model at the end. Stability Selection can be interpreted as meta-algorithm for which algorithms like the usual Lasso but also Lasso variants and Boosting (\cite{hofner15}) can enter as base algorithm for feature selection. 

All non-ensemble non-robust feature selection algorithms, including standalone Boosting models, usually get distorted once the training data are contaminated. Safeguarding against such contamination is done by methods of Robust Statistics (e.g., \cite{maronna}, \cite{huber}, \cite{hampel}, \cite{rieder}). Although outlier detection strategies (e.g., \cite{rocke}, \cite{shieh}, \cite{filz08}, \cite{rous11}, \cite{rous16}) are important for data cleaning so that a classical algorithm may be applied, robust learning algorithms can directly cope with contaminated data. The main idea (see \cite{huber}, \cite{hampel}) is to either bound the effects of contamination by replacing unbounded loss functions like the squared loss by a loss function with a bounded derivative or, even better, by a redescender whose derivative tends to zero for large arguments, or to assign less weight to suspicious instances. Robust techniques have successfully entered sparse feature selection like robust Lasso algorithms (\cite{rosset}, \cite{chen10a}, \cite{chen10b}, \cite{chang17}), robust $L_2$-Boosting (\cite{lutz}) and Sparse Least Trimmed Squares (SLTS) (\cite{alfons13}), among several other variants. 

In practice, especially for high-dimensional data, cell-wise outliers (\cite{alqallaf}) are more realistic than case-wise outliers, so single cells are allowed to be contaminated, independently from whether the other cells in the respective instance are contaminated. For high-dimensional data, this can cause each instance to be contaminated while the fraction of outliers, measured by the relative part of contaminated cells and not of contaminated rows, can still remain very low. There are already some algorithms that can cope with cell-wise outliers like the robust clustering algorithm that trims outlying cells (\cite{garcia21}), cell-wise robust location and scatter matrix estimation algorithms (\cite{agos15}, \cite{leung17}) and several regression approaches (\cite{bottmer}, \cite{leung16}, \cite{filz20b}). There are also notable algorithms for detecting and imputing cell-wise outliers like the DDC from \cite{rous16} or the cellHandler from \cite{rous19b}. 

A common robustness measure is given by the breakdown point. In contrast to the prominent influence curve which quantifies the local robustness of an estimator, i.e., only allowing for an infinitesimal fraction of the data to be contaminated, the breakdown point (BDP), introduced in \cite[Sec. 6]{hampel71} in a functional version and in \cite{huber83} in a finite-sample version, studies the global robustness of an estimator. The finite-sample BDP from \cite{huber83} quantifies the minimum fraction of instances in a data set that guarantees the estimator to ``break down'' when being allowed to be contaminated arbitrarily while the functional BDP essentially quantifies the minimum Prokhorov distance of the ideal and contaminated distribution that leads to such a breakdown. There has already been a lot of work on BDPs, see for example \cite{rous84}, \cite{rous85}, \cite{davies93} and \cite{hubert97}, \cite{genton98}, \cite{becker99}, \cite{gather} or \cite{hubert08} which cover location, scale, regression, spatial and multivariate estimators. Recently, BDP concepts for classification (\cite{zhao18}), multiclass-classification (\cite{qian}) and ranking (\cite{TW21}) have been proposed. Another type of breakdown point which relates the quotient of the number of variables and observations to the sparsity of the true model was studied in \cite{donoho06}, illuminating that in high-dimensional settings, classical model selection procedures like Lasso can only find a reliable model provided that the true model is sufficiently sparse. 

Despite the recent successes of robust model selection methods and path-breaking theoretical results concerning model selection consistency (see, e.g., \cite{bu}), the question how contamination affects variable selection is seldomly addressed, leaving the connection of the paradigms of sparse model selection (select a small fraction of columns), stability (select variables that are appropriate for the majority of the data, i.e., a majority of the rows) and robustness (focus on a certain ``clean'' majority of the rows) still opaque. Although the frequency of contamination in real data is high and although the issue that non-aggregated feature selection models are usually unstable and tend to overfit is well-known, combining both stability and robustness seems to have rarely been considered in literature so far. A robust Stability Selection based on cleaned data has been introduced in \cite{uraibi}. \cite{uraibi19} apply a weighted LAD-Lasso (\cite{arslan}) as base algorithm for the Stability Selection where the weights are computed according to a robust distance in order to downweigh leverage points. \cite{park19} proposed a robust Stability Selection where the term ``robust'' however refers to an immunization against a specific regularization parameter of the underlying Lasso models, therefore considers a different goal than we do. Notable work on aggregating robust estimators has been done in \cite{salibian02}, \cite{salibian08} who propose a linear update rule for robust MM-estimators in order to avoid computing it on each of the drawn Bootstrap samples, see also \cite{salibian06a} for fixed-design regressor matrices and \cite{salibian06b} for robust PCA estimators. While these techniques do not consider variable selection, \cite{salibian08b} extend their method to variable selection where a backward selection strategy according to a minimization of the expected prediction error is applied. 

We aim at making a tiny step towards the connection of robustness, stability and sparsity by introducing the variable selection breakdown point, which describes the number of outliers that can make variable selection completely unreliable, and the Stability Selection breakdown point, which corresponds to a sufficiently high probability that a stable model becomes completely unreliable. We study the relative robustness improvement that a Stability Selection grants, compared to a single model selection algorithm. It turns out that a rank-based Stability Selection where a fixed number of best variables enters the stable model is generally more robust than a threshold-based Stability Selection where all variables whose aggregated selection frequency exceeds the threshold enters the stable model. We propose a trimmed Stability Selection (\textbf{TrimStabSel}) and investigate its performance on a large variety of simulated data in comparison with a non-robust Stability Selection, where $L_2$-Boosting, LogitBoost and SLTS are used model selection algorithms. The numerical results show that even an extremely small cell-wise contamination rate can already have a severe impact on variable selection. Our TrimStabSel can in particular be recommended in settings where the contamination rate is expected to be low but non-zero.

\textbf{Our contribution is threefold:} \textbf{i)} We propose BDP definitions for (stable) variable selection; \textbf{ii)} We propose (oracle) outlier schemes that can completely distort model selection with usually very few outliers; \textbf{iii)} We lift the popular concept of trimming from single instances to whole models that allows for contamination rates exceeding 50\% while maintaining the 50\%-bound for the standard BDPs of the underlying model selection algorithms.

This paper is organized as follows. Section \ref{prelim} compiles relevant notions from Robust Statistics which are contamination models and the breakdown point. We also recapitulate the concept of Stability Selection. Section \ref{bdpsec} is devoted to the definition of our variable selection BDP and to a discussion on resampling BDPs which are required for our Stability Selection. Our BDP concept for Stability Selection is introduced in Section \ref{stabbdpsec}. Section \ref{TrimStabSelsec} presents the Trimmed Stability Selection. Section \ref{simsec} provides a detailed simulation study that compares the performances of the standard Stability Selection and the Trimmed Stability Selection on simulated data with $L_2$-Boosting, LogitBoost and SLTS as model selection algorithms.

\section{Preliminaries} \label{prelim}

Let $D \in \R^{n \times (p+k)}$ be the data matrix consisting of a regressor matrix $X \in \R^{n \times p}$ and a response matrix $Y \in \R^{n \times k}$. If univariate responses are considered, $Y \in \R^n$ is a response column. We denote the $i$-th row of $X$ by $X_i$ and the $j$-th column of $X$ by $X_{\cdot,j}$. The $i$-th row of $Y$ is denoted by $Y_i$ and for $k=1$, $Y_i$ denotes the $i$-th component of $Y$. Let $n_{sub}<n$ always be the number of instances in subsamples resp. Bootstrap samples.  

We start by recapitulating contamination models and the breakdown point concept. The last subsection is devoted to Stability Selection.

\subsection{Contamination models} 

We begin with the definition of contamination balls (\cite[Sec. 4.2]{rieder}). 

\begin{Def} \label{contball} Let $(\Omega, \mathcal{A})$ be a measurable space. Let $\mathcal{P}:=\{P_{\theta} \ | \ \theta \in \Theta\}$ be a parametric model where each $P_{\theta}$ is a distribution on $(\Omega, \mathcal{A})$ and where $\Theta \subset \R^p$ is some parameter space. Let $P_{\theta_0}$ be the ideal distribution (``model distribution'', \cite[Sec. 4.2]{rieder}). Then, a \textbf{contamination model} is the set of all distributions given by the system $\mathcal{U}_*(\theta_0):=\{U_*(\theta_0,r) \ |\ r \in [0,\infty[\}$ of \textbf{contamination balls} $U_*(\theta_0,r)=\{Q \in \mathcal{M}_1(\mathcal{A}) \ | \ d_*(P_{\theta_0},Q) \le r \}$ for the set $\mathcal{M}_1(\mathcal{A})$ of probability distributions on $\mathcal{A}$. The radius $r$ is also called ``contamination radius''. \end{Def}

To emphasize different standard types of contamination balls, a subscript (here represented by the ``$*$'') is usually added to concretize the contamination model.

\begin{ex} \label{convcont} A convex contamination model $\mathcal{U}_c(\theta_0)$ is the system of contamination balls \begin{center} $ \displaystyle U_c(\theta_0,r)=\{(1-r)_+P_{\theta_0}+\min(1,r)Q \ | \ Q \in \mathcal{M}_1(\mathcal{A}) \}. $ \end{center} \end{ex} 

Many other metrics have entered the theory of Robust Statistics, see \cite[Sec. 4.2]{rieder} for an overview. While in standard estimation problems like estimating a location or scale parameter one just has a sample of observations, there is far more flexibility to define outliers in for example regression settings. See \cite{kohl} and \cite{rieder} for more details. 

\begin{rem}[\textbf{Regression outliers}] If we only allow for contaminating the responses while maintaining the regressors, we speak of $Y$-outliers while contamination in the regressors is referred to as $X$-outliers. \end{rem}


In the contamination model in Def. \ref{contball}, we considered case-wise (row-wise/instance-wise) outliers, i.e., either a whole row in the regressor matrix or in the response is contaminated or not. A more realistic scenario where the entries (cells) of the regressor matrix are allowed to be perturbed independently (cell-wise outliers) has been introduced in \cite{alqallaf}. See also \cite{agos15} for the notation. 

\begin{Def} \label{cell-wise} Let $X_i \sim P_{\theta_0}$ for $i=1,...,n$ where $P_{\theta_0}$ is a distribution on some measurable space $(\Omega,\mathcal{A})$ and let $X$ be the $n \times p$-matrix consisting of the $X_i$ as rows. Let $U_1,...,U_p \sim Bin(1,r)$ i.i.d.. Then the \textbf{cell-wise convex contamination model} considers all sets \begin{center} $ \displaystyle U^{cell}(\theta_0,r):=\{Q \ |\ Q=\mathcal{L}(UX+(1-U) \tilde X)\} $ \end{center} where $\tilde X \sim \tilde Q$ for any distribution $\tilde Q$ on $(\Omega,\mathcal{A})$ and the matrix $U$ with diagonal entries $U_i$. \end{Def} 

\begin{rem} \cite{alqallaf} pointed out that if all $U_j$ are perfectly dependent, one either gets the original row or a fully contaminated row which is just the classical convex contamination model in Def. \ref{convcont}. \end{rem}

\begin{rem} Note that for response matrices, one can similarly construct cell-wise outliers. \end{rem}

In the cell-wise contamination model, the curse of dimensionality becomes drastically more apparent since the probability that at least one case is contaminated grows with the dimension $p$. Note that a single contaminated cell already makes an observation an outlier. This has been pointed out in \cite{croux14}, \cite{croux15} who concentrated on precision matrix estimation, see also \cite{loh15b}.

\subsubsection{The breakdown point concept} \label{bdpsec} 

The goal of Robust Statistics is to provide robust estimators, i.e., estimators that tolerate a certain amount of contaminated data without being significantly distorted. As for the term ``robustness'', in this work, we use the global robustness concept that allows for a large fraction of the data points being contaminated arbitrarily. The minimum fraction of outliers that can lead to a breakdown of the estimator is called the breakdown point (BDP) of this particular estimator. The finite-sample BDP of \cite{huber83} is defined as follows. 

\begin{Def} Let $Z_n$ be a sample consisting of instances $(X_1,Y_1),...,(X_n,Y_n)$. The \textbf{finite-sample breakdown point} of the estimator $\hat \beta$ is defined by \begin{equation} \label{fsbdp} \varepsilon^*(\hat \beta,Z_n)=\min\left\{\frac{m}{n} \ \bigg| \ \sup_{Z_n^m}(||\hat \beta(Z_n^m)||)=\infty \right\} \end{equation} where $Z_n^m$ denotes any sample with $(n-m)$ instances in common with the original sample $Z_n$, so one can arbitrarily contaminate $m$ instances of $Z_n$, and where $\hat \beta(Z_n)$ is the estimated coefficient on $Z_n^m$. \end{Def} 

In this definition, one implicitly assumes that $\beta \in \R^p$. If w.l.o.g. a compact set $\Theta \subset \R^p$ is considered, one may define a breakdown in the sense that $\hat \beta$ is located at the boundary of $\Theta$. It has been discussed in \cite{davies05} whether attaining a boundary value should be considered as a breakdown. See the rejoinder of \cite{davies05} for some discussions on that topic where the authors suggest a metric which becomes infinity once one of the boundary values is attained.

There are already a lot of BDP concepts in literature, for example for regression (\cite{stromberg} and \cite{sakata}), for location-scale estimators (\cite{sakata98}), for time series (\cite{genton03}, \cite{genton03b}) and for variogram estimators (\cite{genton98}). \cite{donoho06} propose a BDP for model selection while \cite{kanamori} studied the BDP for SVMs. A BDP variant for clustering was introduced in \cite{hennig08} who called their concept the ``dissolution point''. The situation of heavy-tailed original data has been addressed in \cite{horbenko12} with the expected BDP that takes the ideal distribution of the original data into account. 

\subsection{Stability Selection} \label{stabselsec}

The Stability Selection is an ensemble model selection technique that has been introduced in \cite{bu10}, mainly with the goal to reduce the number of false positives (non-relevant variables that are selected by the algorithm) and also motivated by the fact that the true predictor set $S_0 \subset \{1,...,p\}$ is often not derivable by applying a single model selection procedure. In short, one draws $B$ subsamples from the data of usually around $n/2$ instances and performs a model selection algorithm on each subsample which leads to a set $\hat S^{(b)} \subset \{1,...,p\}$ of selected variables for each $b=1,...,B$. The next step is to aggregate the selection frequencies of all variables, i.e., the binary indicators whether a particular variable has been selected in a specific predictor set. More precisely, one computes $\hat \pi_j:=\frac{1}{B}\sum_{b=1}^B I(j \in \hat S^{(b)})$ for all $j$. 

In the original Stability Selection from \cite{bu10}, one defines a threshold $\pi_{thr}$ based on an inequality derived in \cite[Thm. 1]{bu10} so that the stable set then consists of all variables $j$ for which $\hat \pi_j \ge \pi_{thr}$. There are some variants of this Stability Selection, most notably the one in \cite{hofner15} that makes it applicable for Boosting while the original one is tailored to algorithms that invoke a regularization term like the Lasso or the Graphical Lasso. An excellent implementation of the Stability Selection can be found in the $\mathsf{R}$-packages \texttt{mboost} (\cite{mboost}, \cite{hofner14}, \cite{hothorn10},  \cite{bu07}, \cite{hofner15}, \cite{hothorn06}) and \texttt{stabs} (\cite{stabs}, \cite{hofner15}, \cite{mayr17}). 

As for the selection of the stable set according to the aggregated selection frequencies, another paradigm that defines a number $q$ of variables that have to enter the stable model so that the $q$ variables with the highest selection frequencies are chosen has been suggested in literature (e.g., \cite{zhou13}, \cite{TW22}) since the threshold-based approach is less intuitive for the user due to the number of stable variables not being predictable in the first place.

\section{Breakdown of variable selection} \label{bdpsec}

In this section, we first define a BDP for variable selection. Based on this definition, we discuss why this BDP may be very small and outline the path from the robustness of single algorithms concerning variable selection to ensembles of such algorithms.

\subsection{Variable selection breakdown point} 

\cite{donoho06} already provided a very insightful work where the notion of a ``breakdown of model selection'' has been introduced. They computed phase diagrams that show under which configurations of the dimensionality of the data and the sparsity level of the true underlying model a successful model selection is possible. More precisely, they derive that the underlying model has to be sufficiently sparse, expressed in the fraction of the true dimensionality $q$ of the model and the number $n$ of observations, the question whether model selection is possible depends on the fraction $n/p$ for $p$ being the number of predictors. 

Our idea also considers to compute a breakdown for model selection, but we restrict ourselves to the standard setting of Robust Statistics where we want to examine how many outliers can be tolerated for model selection. In order not to confuse our BDP concept with the one in \cite{donoho06}, we call our concept the variable selection breakdown point (VSBDP). 

\begin{Def} \label{vsbdp} Let $D$ be a data set with $n$ instances and predictor dimension $p$. Let $k \in \N$ be the dimension of the responses. \\
\textbf{a)} The \textbf{case-wise variable selection breakdown point (case-VSBDP)} is given by \begin{equation} \frac{m^*}{n}, \ \ \ m^*=\min\{m \ | \ \hat \beta_j(Z_n^m)=0 \ \forall j: \beta_j \ne 0\} \end{equation} where $Z_n^m$ again denotes any sample that has $(n-m)$ instances in common with $Z_n$. \\ 
\textbf{b)} The \textbf{cell-wise variable selection breakdown point (cell-VSBDP)} is given by \begin{equation} \frac{\tilde m^*}{(p+k)n}, \ \ \ \tilde m^*=\min\{m \ | \ \hat \beta_j(\tilde Z_{cell}^m)=0 \ \forall j: \beta_j \ne 0\} \end{equation} where $\tilde Z_{cell}^m$ denotes the data set where $m$ cells can be modified arbitrarily and where all other cell values remain as in the original data. \end{Def} 

In other words, the VSBDP quantifies the relative number of rows resp. cells that have to be contaminated in order to guarantee that none of the relevant variables are selected. The fraction of outlying cells as a breakdown measure has for example already been considered in \cite{velasco}. Let us now formulate a very simple but important result.

\begin{thm} \label{vsbdpthm} Let $q \le p$ be the true dimension of the underlying model and let again $k$ be the response dimension and $n$ be the number of instances in the data set. \\
\textbf{a)} Then the cell-VSBDP is at most $\frac{\min(q,k)}{p+k}$.\\
\textbf{b)} Let $p=p(n)$ so that it grows when $n$ grows. If $q$ or $k$ stays constant, the asymptotic cell-VSBDP is zero. \\
\textbf{c)} Let $p=p(n)$ and $q=q(n)$ so that both quantities grows when $n$ grows. Let contamination only be allowed on the predictor matrix. Then the asymptotic cell-VSBDP is given by $\lim_{n \rightarrow \infty}\left(\frac{q(n)}{p(n)}\right)$. \end{thm}

\begin{bew} \textbf{a)} If $q< k$, just replace the entries $X_{ij}$ for all $i$ and for all $j$ corresponding to the relevant variables by zeroes, so that one has to modify $qn$ cells of the data set, more precisely only of the predictor matrix. Then the originally relevant columns remain without any predictive power and therefore will not be selected. If $q \ge k$, replace all entries of the response matrix with zeroes which would result in an empty model since no predictor column remains correlated with the response, requiring $kn$ outlying cells. \\
\textbf{b)+c)} Directly follows from a). \end{bew} \vspace{-1cm} \begin{flushright} $_\Box $ \end{flushright}

This is a universal result, regardless of the data structure or the applied algorithms. We want to point out that the classical understanding of robustness would only consider the estimated coefficients themselves and call an estimator robust if the coefficients stay bounded. However, if the non-zero coefficients correspond to non-relevant variables, the learning procedure results in a robust fit on noisy variables which will definitely have poor prediction quality for out-of-sample data. Our analysis is based on an interplay between model selection and coefficient estimation so that the ultimate goal is to achieve both sub-goals in order to get a reliable model. Let us therefore pose the following provoking statement: \textit{From the perspective of retrieving the correct model, all robust regression and classification models are doomed to have a breakdown point less than $q/p$.}

There are methods that solve regression problems by first identifying cell-wise outliers in the data matrix corresponding to a regression problem (\cite{agos15}, \cite{leung16}) and by down-weighting them. In general, the column structure of the contamination scheme corresponds to a perfect correlation of the binomial variables in the cell-wise convex contamination model in Def. \ref{cell-wise}, i.e., all $U_{ij}$, $i=1,...,n$, are perfectly correlated, in contrast to case-wise contamination where the $U_{ij}$, $j=1,...,p$ (or $j=1,...,p+k$ if the response matrix is also considered), are perfectly correlated.\footnote{The  column-wise outliers we defined do not coincide with the ones in \cite{machkour} that also have been termed ``column-wise outliers''. In the reference, the idea is to corrupt columns with column-wise additive outliers.} Note that a zero column would clearly be detected even by hand (although it would not help to reconstruct the original values so that the corresponding variable would not be selected anyway), and other outlier structures like replacing cells by a fixed value could be detected by the outlier detection and imputation procedure (DDC) from \cite{rous16} or the cellHandler from \cite{rous19b}. However, one could try to hide the perturbation by first replacing selected cells of a given column by the mean of the remaining column entries and by adding some random noise to these imputations. 

When having simulated data and a random cell-wise outlier scheme, it is extremely unlikely that such column-wise outliers that make the true model irretrievable would appear, so this situation should at most very rarely occur in a random simulation as for example done in \cite{filz20b} where first a set of instances is selected and for each of these instances, a random fraction of cells are contaminated. \cite{filz20b} identify the column values with measurements made by a specific sensor, so if this sensor is corrupted across a whole study, one indeed would have column-wise outliers. In the simulations in \cite{filz20b}, they restricted themselves to cell-wise contamination in the regressor matrix while \cite{leung16} and \cite{velasco} also allowed for (a small amount of) contamination in the response vector. From a practical perspective, assuming that contamination cannot appear in the responses is a very heavy assumption since if one allows for all $p$ measurement vectors of the variables to be contaminated, it is difficult to argue why of all things the response measurements should be clean.

From a practical perspective, it is important to make an assumption how much information an attacker has. Since Robust Statistics mainly focuses on deriving robust algorithms that can cope with a certain fraction of contamination but not on attacking algorithms that unveil weaknesses of the algorithms by targetedly producing dangerous contamination, this topic has gained far less attraction than in the Deep Learning community. There, one considers adversarial attacks (\cite{goodfellow}), which minimally perturb instances so that an already trained model makes wrong predictions, and poisoning attacks (\cite{jagielski}), which are perturbations that are injected into the data before the model is trained and which therefore are stronger related to Robust Statistics. However, these perturbations are usually bounded by a matrix norm and not by a contamination radius as in convex contamination schemes which nevertheless has also been considered in Robust Statistics (e.g. \cite{burgard}). The success of such attacks depends on the amount of knowledge that the attacker has both on the victim model as well as on the data and the underlying model. Coming back to our column-wise outliers, an attacker that randomly would inject column-wise outliers to the predictor matrix would hardly contaminate the relevant variables by chance unless it is an oracle attacker that is aware of the relevant variables, maybe due to intercepting and analyzing the data first.

\subsection{Resampling and robustness} 

For Bootstrapping, the probability that in a particular resample there is at least one contaminated instance is given by \begin{equation} a^{Boot}_{1/n,n_{sub},B}(\hat \varepsilon):=1-P(Bin(n_{sub},\hat \varepsilon)=0)^B \end{equation} where $\hat \varepsilon=m/n$ indicates the empirical rate of outlying rows in the data set. Similarly, for subsamples, it is given by \begin{equation} a^{Subs}_{1/n,n_{sub},B}(\hat \varepsilon):=1-P(Hyp(n,n-m,n_{sub})=n_{sub})^B \end{equation} where $Hyp(n,n-m,n_{sub})$ is the hypergeometrical distribution describing the number of clean instances (successes) for $(n-m)=(1-\hat \varepsilon)n$ clean instances in the ``urn''. This observation has already been made in literature, more precisely, in \cite{berrendero} for Bootstrap samples and in \cite{camponovo} for subsampling. \cite[Sec. 3.5]{filz20} also derive the necessary number of Bootstrap samples so that the probability of having at least one clean resample is sufficiently large. This paradigm has been extended for the bagged median in \cite{berrendero}. A general condition so that a bagged robust estimator where the estimator on each resample has the BDP $c$ does break down is given by \begin{equation} a^{Boot}_{c,n_{sub},B}(\hat \varepsilon):=1-P(Bin(n_{sub},\hat \varepsilon)<\lceil cn_{sub} \rceil)^B \end{equation} for Bootstrapping and by  \begin{equation} a^{Subs}_{c,n_{sub},B}(\hat \varepsilon):=1-P(Hyp(n,n-m,n_{sub}) > \lfloor (1-c)n_{sub} \rfloor)^B \end{equation} for subsampling. This also consistently covers the case of $c=1/n$ since $n_{sub}/n<1$.  

The last two formulae correspond to non-robust resampling aggregation with robust estimators. \cite{berrendero} also considered a robust aggregation procedure called Bragging which was introduced in an earlier version of \cite{bu12} which just replaces the mean aggregation by a median aggregation of the individual estimators. The condition for a breakdown of the bragged estimator changes to at least $\lfloor (B+1)/2 \rfloor$ resamples being sufficiently contaminated to let these resample-individual estimators break down, see also \cite{berrendero}. Similar results are clearly available if one considers trimmed bagging (\cite{croux07}). See also \cite{salibian02} for such a resampling BDP for Bootstrapped robust quantile estimators. 

We now recapitulate the resampling breakdown point introduced in \cite{berrendero} and use a slightly modified definition to make it more consistent with other BDP concepts.

\begin{Def} Let $n$ be the number of observations in a data set and let $B$ be the number of Bootstrap resp. subsamples with $n_{sub}<n$ observations each. Let $c \in [0,1[$ be the BDP of each estimator applied on the individual Bootstrap samples resp. subsamples. For a tolerance level $\alpha \in [0,1[$, the \textbf{$\alpha$ resampling BDP for Bootstrap resp. subsampling} is given by \begin{equation} (\varepsilon^{\perp}(c,n_{sub},B,\alpha))^*:=\inf\{\varepsilon \in \{0,1/n,...,1\} \ | \ a^{\perp}_{c,n_{sub},B}(\varepsilon)> \alpha \} \end{equation} for $\perp \in \{Boot, Subs\}$. \end{Def}

This resampling BDP defines the maximum fraction of contaminated instances in the data so that the probability that a mean aggregation breaks down exceeds the tolerance level. This definition is important since it lifts the worst-case BDP concept to a probabilistic concept that respects that the worst case is often very unlikely and that solely reporting it would be too pessimistic. 

Now, after this exposition, let us first emphasize that the idea from the rejoinder of \cite{davies05} to map the boundary values of a bounded image set of an estimator to infinite values is extremely important when assessing the effects of resampling on the robustness. \cite{grandval} claim that Bagging never improves the BDP. This is not true if the value set is bounded. We first formally spell out why Bagging is usually not robust, although this fact has already been observed in literature.

\begin{prop} \label{baggingrob} Let $B$ be the number of resamples and let $S_n$ be some estimator, mapping onto an unbounded domain, for a data set with $n$ observations. If the (classical) BDP of the estimator is $c$, so it is for the bagged estimator. \end{prop}

\begin{bew} Since the BDP of the estimator is $c$, manipulating a relative fraction of $c$ instances in one single resample $b$ suffices to let the estimator $S_n^{(b)}$ on this resample to break down, i.e., the norm of the estimated value is fully controlled by the outliers. Then, as the bagged estimator being the empirical mean of all $S_n^{(b)}$, $b=1,..,B$, its norm is also fully controlled. \end{bew} \vspace{-1.25cm} \begin{flushright} $_\Box $ \end{flushright} 

The proof is rather unusual since it assumes that one can targetedly contaminate a selected resample. Usually, the attacker should only have access to the whole training data set. Then, a probabilistic statement in the spirit of the resampling BDP becomes more appropriate.

\begin{prop} Let $B$ be the number of resamples and let $S_n$ be some estimator, mapping onto an unbounded domain, for a data set with $n$ observations. If the (classical) BDP of the estimator is $c$, the resampling BDP of the bagged estimator equals $(\varepsilon^{\perp}(c,n_{sub},B,\alpha))^*$ for $\perp=Boot$ for Bootstrap resp. for $\perp=Subs$ for subsampling.  \end{prop}

\begin{bew} Evidently, if at least one resample is contaminated to an extent such that the estimator breaks down, due to the unbounded domain and the mean aggregation, the bagged estimator breaks down. By definition of the resampling BDP, this happens with a probability of more than $\alpha$ if the fraction of contaminated rows is $(\varepsilon^{\perp}(c,n_{sub},B,\alpha))^*$. \end{bew} \vspace{-1.25cm} \begin{flushright} $_\Box $ \end{flushright} 

\begin{rem} \textbf{a)} These results can be trivially extended to the case of cell-wise contamination. \\
\textbf{b)} Note that if the fraction of outlying instances/cells is smaller than the case-BDP $c$ resp. the cell-BDP $\tilde c$ of the applied estimator, not letting it break down, resampling induces a certain probability, namely $a^{\perp}_{c,n_{sub},B}(\hat \varepsilon)$, that the bagged estimator indeed suffers a breakdown. One has to keep in mind that this probability may become rather high on specific configurations. Consider Bagging with $B=100$ Bootstrap samples, $n=200$, $n_{sub}=100$ and $m=90$ outlying instances and let $c=0.5$. Then, the probability that a single Bootstrap sample is sufficiently clean is $P(Bin(100,0.45)<50) \approx 0.817$, but the probability that the bagged estimator breaks down is $a^{Boot}_{0.5,100,100}(0.45) \approx 1$. In this case, Bagging would completely de-robustify the estimator. This fact is well-known (see e.g. \cite{salibian08}), but we need it for comparison with Stability Selection, i.e., a bagged variable selection which operates on a finite space. \end{rem}

As for bounded domains however, the situation becomes dramatically different. 

\begin{prop} \label{robbound} \textbf{a)} Let $B$ be the number of resamples and let $S_n$ be some estimator, mapping onto a bounded domain, for a data set with $n$ observations from which no one takes any of the boundary values. If the (classical) BDP of the estimator is $c$, the resampling BDP of the bagged estimator is \begin{equation} \inf\{\varepsilon \in \{0,1/n,...,1\} \ | \ P(Bin(B,P(Bin(n_{sub},\varepsilon)\ge \lceil cn_{sub} \rceil))=B)> \alpha \} \end{equation} for Bootstrapping resp. \begin{equation} \inf\{\varepsilon \in \{0,1/n,...,1\} \ | \ P(Bin(B,P(Hyp(n,n-\varepsilon n,n_{sub}) \le \lfloor (1-c)n_{sub} \rfloor))=B)> \alpha \} \end{equation} for subsampling. \\
\textbf{b)} For Bragging, the resampling BDP becomes \begin{equation} \inf\{\varepsilon \in \{0,1/n,...,1\} \ | \ P(Bin(B,P(Bin(n_{sub},\varepsilon)\ge \lceil cn_{sub} \rceil)) \ge \lfloor (B+1)/2 \rfloor)> \alpha \} \end{equation} for Bootstrapping resp. \begin{equation} \inf\{\varepsilon \in \{0,1/n,...,1\} \ | \ P(Bin(B,P(Hyp(n,n-\varepsilon n,n_{sub}) \le \lfloor (1-c)n_{sub} \rfloor)) \ge \lfloor (B+1)/2 \rfloor)> \alpha \} \end{equation} for subsampling.    \end{prop}

\begin{bew} \textbf{a)} If at least one estimator does not break down, it takes a value in the interior of the domain by assumption. Hence, the aggregated estimated value will also lie in the interior of the domain, so there is no breakdown. A breakdown occurs if and only if all estimators break down to the same boundary value (which, by the worst-case perspective of the BDP concept, can be assumed to be possible). \\
\textbf{b)} When taking the median of the estimated values, it suffices that at least the half estimators have broken down to the same boundary value. \end{bew} \vspace{-1cm} \begin{flushright} $_\Box $ \end{flushright} 

This is an extremely strange result since the median aggregation (and any other trimmed aggregation) is less robust than the mean aggregation. This artifact, resulting from a wrong notion of robustness for estimation, shows that a breakdown has to be defined very carefully if the domain is bounded, providing another argument why the concept from the rejoinder of \cite{davies05} to map the boundary values to infinite values is necessary.

\section{Stability Selection and robustness} \label{stabbdpsec}

Important work on stability of feature selection resp. feature ranking has been done in \cite{nogueira16}, \cite{nogueira17}, \cite{nogueira17b}. \cite{nogueira17b} point out that stable feature selection is either represented by a hard subset selection of the candidate variables or by a ranking of the variable  or individual weights which, given some threshold for the ranks resp. the weights, eventually leads to a subset of variables. The cited works propose similarity metrics in order to quantify the stability of feature selection resp. feature ranking. \cite{nogueira17b} consider the stability of feature selection as a robustness measure for the feature preferences. Note that this robustness notion differs from the definition of robustness in the sense of Robust Statistics. 



\subsection{The Stability Selection BDP}

Having the VSBDP and the resampling BDP defined, we are ready to define a BDP for Stability Selection itself.

\begin{Def} \label{stabbdp} Let $n$ be the number of instances in a data set and let $n_{sub}$ be the number of instances in each resample. Let $B$ be the number of resamples for the Stability Selection and let $R$ be the resampling distribution. Let $S^{stab}(\perp,n_{sub},B,Z_n)$ denote the stable set derived from $B$ resamples according to $\perp=Boot$ or $\perp=Subs$ with $n_{sub}$ instances from the data set $Z_n$. \\
\textbf{a)} Then the \textbf{case-wise Stability Selection BDP (case-Stab-BDP)} for tolerance level $\alpha$ is given by \begin{equation} \varepsilon^*_{Stab}(\perp,c,n_{sub},B,\alpha):=\min\left\{\frac{m}{n} \ \bigg| \ P_R(\forall j \in S_0: j \notin S^{stab}(\perp,n_{sub},B,Z_n^m)) \ge \alpha\right\}  \end{equation} where $c$ represents the case-BDP of the underlying model selection algorithm. \\
\textbf{b)} Similarly, the \textbf{cell-wise Stability Selection BDP (cell-Stab-BDP)} for tolerance level $\alpha$ is given by \begin{equation} \tilde \varepsilon^*_{Stab}(\perp,\tilde c,n_{sub},B,\alpha):=\min\left\{\frac{\tilde m}{n(p+k)} \ \bigg| \ P_R(\forall j \in S_0: j \notin S^{stab}(\perp,n_{sub},B,\tilde Z_{cell}^m)) \ge \alpha\right\} \end{equation} with the cell-BDP $\tilde c$ of the underlying model selection algorithm. \end{Def}

Intuitively, the Stab-BDP denotes the minimum fraction of outliers required so that the probability that the reported stable model does not contain any relevant variable is sufficiently large. As Stability Selection does not aggregate coefficients but just indicator functions, the influence of a model computed on a single resample is bounded by $1/B$ in the sense that the aggregated selection frequencies computed on the original data can be at most distorted by $1/B$ (for certain $j$) if one single resample is (sufficiently) contaminated.

\begin{ex}[Random Lasso] Before we analyze the Stab-BDP, we want to emphasize that this robustification effect gets lost if the model aggregation is based on values from unbounded sets. An exemplar for this model type is the Random Lasso from \cite{wang11} which, for $b=1,...,B$, does not only draw resamples from the rows but also randomly selects several columns from a uniform distribution over the columns. Then, a Lasso model is computed on each resample, leading to coefficients $\hat \beta^{(b)}$. The importance of the columns is then computed by $|\sum_b \hat \beta^{(b)}|/B$ and in a second stage, one draws again resamples from the rows but for each resample one draws columns according to the importance values, computes again the Lasso on each subsample and aggregates the coefficients. Neither the randomization over the columns nor the computation of the importance provides any robustness advantage. One just needs to be able to produce an arbitrarily large coefficient for a non-relevant variable on one single resample. This will cause the importance of this column to be arbitrarily close to 1 in this component and therefore arbitrarily close to 0 in all other components, especially those corresponding to the relevant variables which therefore will not be selected in the second stage. Therefore, the variable selection robustness of Random Lasso can be quantified by the resampling BDP $(\varepsilon^{\perp}(1/n,n_{sub},B,\alpha))^*$ for case-wise outliers and by $(\varepsilon^{\perp}(1/(n(p+1)),n_{sub},B,\alpha))^*$ for cell-wise outliers, provided that one outlying instance resp. cell can \textbf{promote a non-relevant variable}. \end{ex}

\subsection{Can (non-relevant) variables be targetedly promoted in practice?}

Research in this direction has recently been done in the context of poisoning attacks, i.e., injecting contamination into the training data so that one can targetedly manipulate the trained model. In contrast to the BDP concept where one can modify a certain fraction of the data arbitrarily, poisoning attacks generally allow to contaminate all data instances but with a bound on some $||\cdot||_s$-norm for the contamination vector resp. matrix. \cite{lai20f}, \cite{lai20g} prove that poisoning attacks can indeed targetedly affect model selection in the sense that variables selected by the attacker can be suppressed resp. promoted, i.e., the attacker can inject variables to the selected model resp. discard variables. According to their attack scheme, there is essentially no limitation in the number of promoted resp. suppressed variables, but note that they propose $l_s$-attacks for $s \ge 1$, i.e., one does not change a single instance but in principle each instance/cell. Nevertheless, applying their method with an $l_0$-constraint as already done for crafting sparse adversarial attacks (e.g., \cite{su19}, \cite{carlini}) would indeed exactly represent our setting. 

To the best of our knowledge, a concise statement on the impact of such a fraction of outliers to model selection itself like how many relevant variables can be suppressed or how many non-relevant variables can be promoted has not yet been proposed in literature. Therefore, we propose an optimistic and a pessimistic scenario concerning this impact. 

\textbf{Case-wise scenarios:} In this scenario, we assume that for a variable selection method with instance-BDP $c$, a number of $\lceil cn \rceil$ outlying rows in the regressor matrix, the response matrix or the regressor matrix reduced to the relevant columns \begin{itemize} \item is able to promote any subset of non-relevant variables resp. suppress any subset of relevant variables which means that we assume that this outlier fraction can indeed, at least theoretically, cause all relevant variables to be ignored (\textbf{pessimictic case-wise scenario});
\item is able to targetedly suppress all relevant but promote only one single variable (\textbf{optimistic case-wise scenario}). \end{itemize}

\textbf{Cell-wise scenarios:} We assume that for a cell-BDP of $\tilde c$, a fraction of at least $\tilde c$ of outlying cells in the regressor matrix, the response matrix or the regressor matrix reduced to the relevant columns \begin{itemize} \item can promote $\min(p-s_0,\tilde c)$ non-relevant variables where $s_0=|S_0|$ resp. suppress any subset of relevant variables  (\textbf{pessimictic cell-wise scenario}). This scenario is indeed extremely pessimistic since even in \cite{lai20f}, \cite{lai20g}, it seems that one at least has to manipulate the regressor matrix and the response vector which would at least require two outlying cells for an estimator with $\tilde c=0$;
\item is able to targetedly suppress all relevant variables but to promote only one single variable (\textbf{optimistic cell-wise scenario}). \end{itemize}

\subsection{Threshold-based Stability Selection} \label{thrss} 

We now quantify the robustness of threshold-based Stability Selection. We assume that there is a fixed selection of instances resp. cells in the data matrix which are contaminated, so we quantify the probability that the Stability Selection breaks down. 

\begin{thm} \label{thrstabselrob} Let $\pi$ be the threshold for the Stability Selection based on $B$ resamples of size $n_{sub}$ from a data set with $n$ instances, a fixed selection of $m$ of them contaminated, and let $c$ be the BDP of the applied model selection algorithm. Let $\hat \pi_j^+$ be the aggregated selection frequencies on the original data where $\hat \pi_j^+ \ge \pi \ \exists j$ for $j \in S_0$. Then, in both the optimistic and the pessimistic case-wise scenario, the probability that the Stability Selection breaks down is \\
\textbf{i)} given by \begin{equation} \label{hyp1} P(Bin(B,P(Hyp(n,n-m,n_{sub}) \le \lfloor (1-c)n_{sub} \rfloor))>\lceil B(\max(\hat \pi_j^+)-\pi) \rceil)  \end{equation} if the resamples are drawn by subsampling; \\
\textbf{ii)} given by \begin{equation} \label{boot1} P(Bin(B,P(Bin(n_{sub},m/n) \ge \lceil cn_{sub} \rceil))>\lceil B(\max(\hat \pi_j^+)-\pi) \rceil) \end{equation} if the resamples are drawn by Bootstrapping. \end{thm}

\begin{bew} First, note that the selection frequencies of the non-relevant variables are not important here, so one does not distinguish between the pessimistic and the optimistic scenario. According to our assumption, there is no immediate breakdown on the original data. A breakdown is achieved once each relevant variable no longer appears in the stable set, i.e., if all aggregated selection frequencies are below the threshold $\pi$. Regarding the variable $j^*=\argmax_j(\hat \pi_j^+)$, there are more than $\lceil B(\hat \pi_{j^*}^+-\pi) \rceil$ sufficiently contaminated resamples required since each such resample can decrease the aggregated selection frequency by at most $1/B$. Sufficiently contaminated means that at least $\lceil cn_{sub} \rceil$ instances in the resample have to be contaminated. Putting everything together, we get the stated probabilities. \end{bew}  \vspace{-1cm} \begin{flushright} $_\Box $ \end{flushright} 

Note that the column-wise outliers do not have to be considered in the theorem above since it is impossible to have a column-wise contamination rate higher than $c$ without having a case-wise contamination rate higher than $c$ which already is a breakdown, so one does not have to distinguish between contamination in the relevant or non-relevant or response columns here. The situation becomes somewhat different when considering cell-wise contamination as shown in the following theorem where we assume univariate responses.

\begin{thm} \label{thrstabselrobcor} Let $\pi$ be the threshold for the Stability Selection based on $B$ resamples of size $n_{sub}$ from a data set with $n$ instances and let $\tilde c$ be the cell-BDP of the applied model selection algorithm. Let, for $i=1,...,n$, a fixed selection of $c_i$ cells in instance $i$ be contaminated, let $Z_l$, $l=0,1,...,p+1$, be the number of instances for with $c_i=l$ and let $\tilde m:=\sum_i c_i$. Let further $Z_{l'}^{rel}$, $l'=1,...,s_0$, denote the number of instances with $l'$ outlying cells in the relevant columns. Let $m$ be the number of outliers in the response column. Let $\hat \pi_j^+$ be the aggregated selection frequencies on the original data where $\hat \pi_j^+ \ge \pi \ \exists j$ for $j \in S_0$. Then, in both the optimistic and the pessimistic cell-wise scenario, the probability that the Stability Selection breaks down is  \\
\textbf{i)} $1$ if the fraction of cell-wise outliers exceeds $\tilde c$ in the relevant columns or in the response column; \\
\textbf{ii)} $0$ if $c_i \le \lfloor \tilde c(p+1) \rfloor \ \forall i$ resp. $1$ if $c_i>\lfloor \tilde c(p+1) \rfloor \ \forall i$;\\
\textbf{iii)} given by $\min(P_1,P_2,P_3)$ for \begin{center} $ \displaystyle P_1:=P(Bin(B,p_1)>\lceil B(\max(\hat \pi_j^+)-\pi) \rceil) $ \end{center} where \begin{center} $ \displaystyle p_1:=\sum_{z_0,...,z_{p+1}: \sum_l lz_l \ge \lceil \tilde cn_{sub}(p+1) \rceil} f_{Z,n_{sub}}(z_0,...,z_{p+1}) $ \end{center} where $f_{Z,n_{sub}}$ for $Z=(Z_0,...,Z_{p+1})$ represents the probability function of a multivariate hypergeometric distribution with values in $\{(z_0,...,z_{p+1}) \ | \ z_0+...+z_{p+1}=n_{sub}\}$, i.e.,  $f_{Z,n_{sub}}(z_0,...,z_{p+1})$ is the probability that when sampling $n_{sub}$ instances without replacement, one gets $z_0$ out of the $Z_0$ instances without cell-wise outliers, $z_1$ out of the $Z_1$ instances with one cell-wise outlier and so forth; for \begin{center} $ \displaystyle P_2:=P(Bin(B,p_2)>\lceil B(\max(\hat \pi_j^+)-\pi) \rceil) $ \end{center} where \begin{center} $ \displaystyle p_2:=\sum_{\tilde z_0,...,\tilde z_{s_0}: \sum_l l\tilde z_l \ge \lceil \tilde cn_{sub}s_0 \rceil} f_{Z^{rel},n_{sub}}(\tilde z_0,...,\tilde z_{s_0}) $ \end{center} for $Z^{rel}=(Z_0^{rel},...,Z_{s_0}^{rel})$; and for $P_3$ being the quantity in Eq. \ref{hyp1}, if the resamples are drawn by subsampling and if none the conditions in i) or ii) are satisfied; \\
\textbf{iv)} given by $\min(\check P_1,\check P_2,\check P_3)$ for \begin{center} $ \displaystyle \check P_1:=P(Bin(B,\check p_1 >\lceil B(\max(\hat \pi_j^+)-\pi) \rceil) $ \end{center} for \begin{center} $ \displaystyle \check p_1:=\sum_{z_0,...,z_{p+1}: \sum_l lz_l \ge \lceil \tilde cn_{sub}(p+1) \rceil} f_{Mult}(z_0,...,z_{p+1}) $ \end{center} where $f_{Mult}$ represents the density of a multinomial distribution with parameters $n_{sub}$ and $(Z_0/\sum_l Z_l,...,Z_{p+1}/\sum_l Z_l)$; for \begin{center} $ \displaystyle \check P_2:=P(Bin(B,\check p_2 >\lceil B(\max(\hat \pi_j^+)-\pi) \rceil) $ \end{center} for \begin{center} $ \displaystyle \check p_2:=\sum_{\tilde z_0,...,\tilde z_{s_0}: \sum_l l\tilde z_l \ge \lceil \tilde cn_{sub}s_0 \rceil} f_{Mult}(\tilde z_0,...,\tilde z_{s_0}); $ \end{center} and for $\check P_3$ as in Eq. \ref{boot1}, if the resamples are drawn by Bootstrapping and if none the conditions in i) or ii) are satisfied. \end{thm}

\begin{bew} \textbf{i)+ii)} Obvious and already partially discussed earlier. \\
\textbf{iii)} There are three ways how to achieve a breakdown as already mentioned when defining the cell-wise scenarios. $P_1$ quantifies the probability that due to subsampling, the fraction of cell-wise outliers in the whole data matrix becomes at least $\tilde c$ while $P_2$ quantifies the analogous probability for the set of relevant columns. The stated probabilities $p_1$ and $p_2$ quantify that this happens for one resample, so $P_1$ and $P_2$ quantify that this happens in sufficiently many resamples. $P_3$ quantifies the probability that a fraction of at least $\tilde c$ of the responses are contaminated which also can cause the wrong variables to be selected according to the assumption in the cell-wise scenarios. Assuming that $m$ outliers are apparent in the response column, we get the same formula as in Thm. \ref{thrstabselrob}.  \\
\textbf{iv)} As iii). \end{bew}  \vspace{-1.25cm} \begin{flushright} $_\Box $ \end{flushright} 

The Stab-BDP for tolerance $\alpha$ then is the maximum fraction of outliers, i.e., the maximum fraction $m/n$ resp. $\tilde m/(n(p+1))$ such that the respective probability that the Stability Selection breaks down exceeds $\alpha$. 

\begin{rem} \textbf{i)} The cases corresponding to $P_1$ and $P_2$ resp. to $\check P_1$ and $\check P_2$ in Thm. \ref{thrstabselrobcor} have to be considered separately since having only contamination in the relevant columns makes it impossible that the whole predictor matrix is contaminated too much provided $s_0/p<\tilde c$. Similarly, if contamination only occurs in the non-relevant columns, a breakdown can still be possible due to promoting effects of the contamination. \\
\textbf{ii)} In Thm. \ref{thrstabselrobcor}, we considered the case of univariate responses. For multivariate responses with $k$ response columns, one has to distinguish between the seemingly unrelated regression case (\cite{zellner}) where one fits a model for each response column separately and the general case that the response columns are correlated so that the entire response matrix enters as input for a unified model. In the second case, the probability of a breakdown would be $1$ if the relative part of contaminated cells in the response matrix is at least $\tilde c$, in the first case however, the relative part of outliers has to be larger than $\tilde c$ for one resp. for each response column if a breakdown of the set of the resulting $k$ models is defined in the sense that at least one resp. each of the individual models break down. \end{rem}

\subsection{Rank-based Stability Selection} \label{rankss} 

The robustness of the rank-based Stability Selection additionally depends on the aggregated selection frequencies of the non-relevant variables.

\begin{thm} \label{rankstabselrob} Let $q$ be the pre-scribed number of stable variables for the Stability Selection based on $B$ resamples of size $n_{sub}$ from a data set with $n$ instances, $m$ of them contaminated, and let $c$ be the BDP of the applied model selection algorithm. Let $\hat \pi_j^+$ for $j \in S_0$ resp. $\hat \pi_k^-$ for $k \in \{1,...,p\} \setminus S_0$ be the aggregated selection frequencies on the original data where we assume that $\sum_k I(\hat \pi_k^- \ge \hat \pi_j^+)<q-1 \ \exists j \in S_0$. \\
\textbf{a)} In the pessimistic scenario, let $s$ be the number of relevant variables in the top-$q$ variables. Then, w.l.o.g., let $\hat \pi_1^+,...,\hat \pi_s^+$ be the corresponding aggregated selection frequencies of these variables on the original data, similarly, let $\hat \pi_1^-,...,\hat \pi_{q-s}^-$ be the aggregated selection frequencies of the $(q-s)$ non-relevant variables among the top$-q$ variables. Let $\hat \pi_{q-s+1}^-,...,\hat \pi_q^-$ be the aggregated selection frequencies of the next best $s$ non-relevant variables. Then the probability that the Stability Selection breaks down is \\
\textbf{i)} given by \begin{equation} \label{hyp2}  P(Bin(B,P(Hyp(n,n-m,n_{sub}) \le \lfloor (1-c)n_{sub} \rfloor)) > \lceil 0.5B(\max_{j=1,...,s}(\hat \pi_j^+)-\min_{k=q-s+1,...,q}(\hat \pi_k^-)) \rceil)  \end{equation} if the resamples are drawn by subsampling; \\
\textbf{ii)} given by \begin{equation} \label{boot2} P(Bin(B,P(Bin(n_{sub},m/n) \ge \lceil cn_{sub} \rceil)) > \lceil 0.5B(\max_{j=1,...,s}(\hat \pi_j^+)-\min_{k=q-s+1,...,q}(\hat \pi_k^-)) \rceil)  \end{equation} if the resamples are drawn by Bootstrapping.\\
\textbf{b)} In the optimistic scenario, the respective probability \\
\textbf{i)} lies in the interval \begin{equation} \label{hyp3} \begin{split} [P(Bin(B,P(Hyp(n,n-m,n_{sub}) \le \lfloor (1-c)n_{sub} \rfloor)) > \lceil B(\max_{j=1,...,s}(\hat \pi_j^+)-\min_{k=q-s+1,...,q}(\hat \pi_k^-)) \rceil), \\ P(Bin(B,P(Hyp(n,n-m,n_{sub}) \le \lfloor (1-c)n_{sub} \rfloor)) > \lceil 0.5B(\max_{j=1,...,s}(\hat \pi_j^+)-\min_{k=q-s+1,...,q}(\hat \pi_k^-)) \rceil)]  \end{split} \end{equation} if the resamples are drawn by subsampling; \\
\textbf{ii)} lies in the interval \begin{equation} \label{boot3} \begin{split} [ P(Bin(B,P(Bin(n_{sub},m/n) \ge \lceil cn_{sub} \rceil)) > \lceil B(\max_{j=1,...,s}(\hat \pi_j^+)-\min_{k=q-s+1,...,q}(\hat \pi_k^-)) \rceil), \\ P(Bin(B,P(Bin(n_{sub},m/n) \ge \lceil cn_{sub} \rceil)) >\lceil 0.5B(\max_{j=1,...,s}(\hat \pi_j^+)-\min_{k=q-s+1,...,q}(\hat \pi_k^-)) \rceil)] \end{split} \end{equation} if the resamples are drawn by Bootstrapping.\end{thm}

\begin{bew} On the original data, there is no immediate breakdown according to the assumption. A breakdown is achieved once each relevant instance has an aggregated selection frequency smaller than the aggregated selection frequencies of $q$ non-relevant variables.  \\
\textbf{a)} Therefore, it suffices to have more than $\lceil 0.5B(\max_{j=1,...,s}(\hat \pi_j^+)-\min_{k=q-s+1,...,q}(\hat \pi_k^-)) \rceil)$ contaminated resamples since in each one, both the non-relevant variables corresponding to $\hat \pi_k^-$, $k=1,...,q$, can be promoted and at the same time, the relevant variables can be suppressed, so the quantities $\max_{j=1,...,s}(\hat \pi_j^+)$ and $\min_{k=q-s+1,...,q}(\hat \pi_k^-)$ move towards each other by $2/B$ steps per contaminated resample. Hence, after $\lceil 0.5B(\max_{j=1,...,s}(\hat \pi_j^+)-\min_{k=q-s+1,...,q}(\hat \pi_k^-)) \rceil$ such steps, they are equal or have already crossed, so having one more contaminated resample, the formerly best relevant variable definitely has a lower aggregated selection frequency than the formerly $q$-th best non-relevant variable. The probabilities are then computed as in Thm. \ref{thrstabselrob}. \\
\textbf{b)} In the optimistic scenario, we can only targetedly promote one non-relevant variable. Therefore, it depends on the terms $s$, $q$, $\max_j(\hat \pi_j^+)$ and all $\hat \pi_k^-$ for $k=1,...,q$, so one cannot make a universal precise statement. However, there are two extreme cases. If for \begin{center} $ \displaystyle\max_{j=1,...,s}(\hat \pi_j^+)-\min_{k=q-s+1,...,q}(\hat \pi_k^-)=:\Delta $ \end{center} the difference between $\max_j(\hat \pi_j^+)$ and $\hat \pi_k^-$ for $(q-1)$ indices $k$ from $\{1,...,q\}$ is exactly $\Delta/2$, except for $k^*:=\argmin_{k=q-s+1,...,q}(\hat \pi_k^-)$ for which it is $\Delta$, then more than $0.5 \lceil B\Delta \rceil$ contaminated samples suffice for a breakdown if the variable corresponding to $\hat \pi_{k^*}^-$ is promoted in each of these resamples since the same reasoning as in a) applies, i.e., the selection probabilities for the worst non-relevant and best relevant variable move towards each other with a step size of $2/B$. The other extreme case is that all $\hat \pi_k^-$ are equal. Then, if $s>\lceil B\Delta \rceil$, even after promoting each of these non-relevant variables in one single contaminated resample does not suffice for a breakdown since there will be at least one remaining one whose aggregated selection frequency was still not promoted. In that case, we can treat $\hat \pi_k^-$ (in general, $\min_{k=q-s+1,...,q}(\hat \pi_k^-)$) as threshold so that the results from Thm. \ref{thrstabselrob} are applicable, so the relevant variables have to be suppressed in so many resamples such that the selection frequency of the best of them finally crosses the threshold. \end{bew}  \vspace{-0.5cm} \begin{flushright} $_\Box $ \end{flushright}

\begin{cor} \label{rankstabselrobcor} Let $q$ be the pre-scribed number of stable variables for the Stability Selection based on $B$ resamples of size $n_{sub}$ from a data set with $n$ instances and let $\tilde c$ be the cell-BDP of the applied model selection algorithm. Let, for $i=1,...,n$, $c_i$ cells in instance $i$ be contaminated, let $Z_l$, $l=0,1,...,p+1$, be the number of instances for with $c_i=l$ and let $\tilde m:=\sum_i c_i$. Let further $Z_{l'}^{rel}$, $l'=1,...,s_0$, denote the number of instances with $l'$ outlying cells in the relevant columns. Let $\hat \pi_j^+$ for $j \in S_0$ resp. $\hat \pi_k^-$ for $k \in \{1,...,p\} \setminus S_0$ be the aggregated selection frequencies on the original data where we assume that $\sum_k I(\hat \pi_k^- \ge \hat \pi_j^+)<q-1 \ \exists j \in S_0$. Then, the probability that the Stability Selection breaks down is \\
\textbf{a)} $1$ if the fraction of cell-wise outliers exceeds $\tilde c$ in the relevant columns or in the response column; \\
\textbf{b)} $0$ if $c_i \le \lfloor \tilde c(p+1) \rfloor \ \forall i$ resp. $1$ if $c_i>\lfloor \tilde c(p+1) \rfloor \ \forall i$;\\
\textbf{c)} Let none of the conditions in a) or b) hold. In the pessimistic scenario, let the notation and assumptions be as in part a) of Thm. \ref{rankstabselrob}. Then, the probability that the Stability Selection breaks down is \\
\textbf{i)} given by $\min(P_1,P_2,P_3)$ for \begin{center} $ \displaystyle P_v:=P(Bin(B,p_v) > \lceil 0.5B(\max_{j=1,...,s}(\hat \pi_j^+)-\min_{k=q-s+1,...,q}(\hat \pi_k^-)) \rceil) $ \end{center} for $v=1,2$ with $p_1$ resp. $p_2$ as in Thm. \ref{thrstabselrobcor} and for $P_3$ as in Eq. \ref{hyp2} if the resamples are drawn by subsampling; \\
\textbf{ii)} given by $\min(\check P_1,\check P_2,\check P_3)$ for \begin{center} $ \displaystyle \check P_v:=P(Bin(B,\check p_v) > \lceil 0.5B(\max_{j=1,...,s}(\hat \pi_j^+)-\min_{k=q-s+1,...,q}(\hat \pi_k^-)) \rceil) $ \end{center} with $\check p_1$ and $\check p_2$ from Thm. \ref{thrstabselrobcor} and $\check P_3$ as in Eq. \ref{boot2} if the resamples are drawn by Bootstrapping.\\
\textbf{d)} Let none of the conditions in a) or b) hold. In the optimistic scenario, the respective probability \\
\textbf{i)} lies in one of the intervals \begin{center} $ \displaystyle [P(Bin(B,p_v) > \lceil B(\max_{j=1,...,s}(\hat \pi_j^+)-\min_{k=q-s+1,...,q}(\hat \pi_k^-)) \rceil), P(Bin(B,p_v) \ge \lceil cn_{sub} \rceil)) > \lceil 0.5B(\max_{j=1,...,s}(\hat \pi_j^+)-\min_{k=q-s+1,...,q}(\hat \pi_k^-)) \rceil)] $ \end{center} for $v=1,2$ with $p_1$ and $p_2$ from Thm. \ref{thrstabselrobcor} or in the interval given in Eq. \ref{hyp3} if the resamples are drawn by subsampling; \\
\textbf{ii)} lies in one of the intervals \begin{center} $ \displaystyle [ P(Bin(B,\check p_v) > \lceil B(\max_{j=1,...,s}(\hat \pi_j^+)-\min_{k=q-s+1,...,q}(\hat \pi_k^-)) \rceil), P(Bin(B, \check p_v) > \lceil 0.5B(\max_{j=1,...,s}(\hat \pi_j^+)-\min_{k=q-s+1,...,q}(\hat \pi_k^-)) \rceil)] $ \end{center} for $v=1,2$ with $\check p_1$ and $\check p_2$ from Thm. \ref{thrstabselrobcor} or in the interval given in Eq. \ref{boot3} if the resamples are drawn by Bootstrapping.\end{cor}

\begin{bew} Parts a), b) and c) directly follow from Thm.s \ref{rankstabselrob} and \ref{thrstabselrobcor}. Statement d) is more tricky since one has intervals due to the argumentation in part b) of the proof of Thm. \ref{rankstabselrob}. For a concrete data set and a concrete model selection algorithm, one value in the respective intervals is realized so that the probability of a breakdown of the Stability Selection is the minimum.  \end{bew}  \vspace{-1cm} \begin{flushright} $_\Box $ \end{flushright}

\subsection{Implications of the theorems}

We learned from the results in the previous subsections that Stability Selection does not suffer from the numerical instabilities (cf. \cite{salibian08} for this notion) of resampling that can lead to some resamples having a larger fraction of outlying instances/cells than the whole training data to that extent as standard bagged estimators do.

Apart from the simple observation that robustness is increased by Stability Selection, we can extract one interesting recommendation for the Stability Selection variant. For the threshold-based Stability Selection, the robustness depends on the difference of the aggregated selection frequency of the best relevant variable and the threshold. It is however not evident that there even exists such a variable. This problem has been studied for example in \cite{TW22} for noisy data. Due to the combination of a high noise level, a large number $p$ of candidate variables and a rather sparse true model, one often faces the situation that no variable (including non-relevant variables) can pass the threshold, leading to an empty model, i.e., an immediate breakdown. The rank-based variant circumvents this issue. Moreover, when considering contamination, achieving a breakdown does not only require to let the aggregated selection frequencies drop below the threshold, but their ranks have to drop below $q$ if the best $q$ variables enter the stable model, which is more difficult. Therefore, we recommend to prefer the rank-based Stability Selection over the threshold-based Stability Selection in regard of variable selection robustness. Note that this does not contradict \cite[Thm. 1]{bu10} as the number $q$ can be fixed empirically so that an appropriate bound is achieved, where the aggregated selection probability of the $q$-th best variable replaces the universal threshold. 

In contrast to a simple bagged estimator, Stability Selection allows for more than the half of the instances/cells being contaminated without violating equivariance properties (cf. \cite{davies05}) of the underlying algorithm that prevent BDPs from exceeding 0.5. For example, in machine learning, especially regression, one assumes that there is an underlying model from which the clean data have been generated. Even if the relative fraction of outlying instances exceeds a half, say it is 60\%, then it is nevertheless a desirable goal to infer the model which describes the correspondence structure of the responses and the predictors of the clean observations. There is no qualitative hindrance to aim at finding the underlying model due to the standard assumption that outliers \textbf{do not have structure} and just stem from some unknown, arbitrary distribution. 

In order to clarify the argumentation, consider the awkward case-wise contamination situation that the outliers had structure, i.e., additionally to the underlying model $f: \mathcal{X} \rightarrow \mathcal{Y}$ that relates the responses and the predictors of the clean observations one had another model $g: \mathcal{X} \rightarrow \mathcal{Y}$ from which the outlying instances are generated. In this artificial setting, one would indeed try to infer $g$ instead of $f$ and treat the actual clean instances as the outliers, more precisely, as the outliers w.r.t. the model $g$. 

The main problem, even for a Stability Selection, would be that the probability to draw resamples that are sufficiently contaminated to let the applied algorithm break down would become rather high if the fraction of outlying instances already exceeds 0.5 on the original data. Therefore, an additional robustification step is necessary which will be discussed in the next section.

\section{Robustifying Stability Selection} \label{TrimStabSelsec}

Let us recapitulate the SLTS algorithm from \cite{alfons13}: It computes the loss function $L$ (the squared loss) for all instances and defines the ``clean subset'' that enters the next iteration as the $h<n$ instances with the lowest in-sample losses. By iterating this strategy and by immunizing against running into a local minimum of the trimmed objective $\sum_{i=1}^h (L(Y,X,\beta))_{i:n}$ by starting with multiple initial configurations, the SLTS outputs essentially the optimal sparse coefficients according to the empirical trimmed risk.

Regarding Stability Selection, the instances in a resample cannot be treated individually since one aims at aggregating column information from whole resamples. One can however individualize the resamples themselves and the corresponding selected sets of variables. In other words, we aim at lifting the trimming concept from instances in estimation problems to resamples in a model aggregation problem.

\subsection{Related approaches} 

A more flexible approach than Bragging (\cite{bu12b}) is trimmed Bagging (\cite{croux07}) where the trimmed mean of the estimators instead of their median is computed. In contrast to \cite{croux07}, we do not bag classifiers nor other learners but predictor sets themselves.  A close approach to our intended Trimmed Stability Selection is the approach in \cite{zhang19d} where general ensemble variable selection techniques are pruned. Their motivation came from the perspective of strongly correlated members of the ensemble in contrast to the perspective of robustness as in our work. They extend the approach in \cite{zhang17} where the members with the highest prediction errors are trimmed in the sense that \cite{zhang19d} define a so-called variable selection loss which is a squared distance between the aggregated selection frequencies (called ``importances'' there) and a reference importance vector. \cite{zhang19d} are aware of the fact that the true importance vector is not known so they approximate it empirically.

A very important related algorithm is the fast and robust Bootstrap (FRB), initially introduced in \cite{salibian02}, \cite{salibian08} for regression MM-estimators. The main idea is that standard Bootstrapping of robust estimators would take a long computational time and that a simple uniform Bootstrapping may would result in having resamples with more than half of the instances being contaminated with a considerable probability. They derived a linear update rule for the robust estimator that downweights outlying instances. These weights are derived by a robust MM-estimator which allows for identifying such instances according to the absolute value of their regression residual. Model selection using the FRB has been proposed by \cite{salibian08b}. The idea is to approximate the prediction error of the model built on a particular set of predictors by using FRB and to select the predictor set for which the prediction error was minimal. Their strategy suffers from the lack of scalability for data sets with a large number of predictors since they have to compute the prediction error for all possible models whose total number is $2^p-1$, and even their backward strategy become infeasible for high $p$. 

\subsection{Trimmed Stability Selection}

In this paper, we do not consider validation data since assuming that the training data are contaminated while the validation data are clean would be awkward, so we intend to measure the quality of the resample-specific models on an in-sample loss basis, similarly as outlying instances are detected by their individual in-sample loss as for example in \cite{alfons13}. More precisely, if the contamination of a certain resample has caused the algorithm to select wrong variables, the in-sample loss should be high compared to another resample where a sufficiently well model has been selected which contains enough of the true predictors. 

Our trimmed Stability Selection works as follows. We generate $B$ resamples of size $n_{sub}$ from the training data and apply the model selection algorithm, for example, Lasso or Boosting, on each resample which selects a model $\hat S^{(b)} \subset \{1,...,p\}$ for $b=1,...,B$ and which computes coefficients $\hat \beta^{(b)}$. Let $I^{(b)}$ be the index set of the rows that have been selected by the resampling algorithm, i.e., $I^{(b)} \in \{1,...,n\}^{n_{sub}}$ for Bootstrapping resp. $I^{(b)} \subset \{1,...,n\}$ with $|I^{(b)}|=n_{sub}$ and $I_k^{(b)} \ne I_l^{(b)} \ \forall k \ne l, k,l=1,...,n_{sub}$, for subsampling. For the $b$-th resample, we compute the in-sample loss \begin{equation} \label{insample} L^{(b)}=\frac{1}{n_{sub}} \sum_{i \in I^{(b)}} L(Y_i,X_i\hat \beta^{(b)})  \end{equation} for a loss function $L: \mathcal{Y} \times \mathcal{Y} \rightarrow [0,\infty[$. Although, in regression, losses are usually continuous, we cannot exclude that there are ties, for example, if two resamples are identical or if a trimmed loss is used. In this case, we perform random tie breaking. Note that there is no evidence that the models computed on identical resamples are identical since the model selection algorithm may include stochastic components like the SLTS whose result depends on the randomly chosen initial configurations. 

For the set \begin{equation} I^{trim}(\gamma):=\left\{b \in \{1,...,B\} \ \bigg| \ \sum_{b'=1}^B I(L^{(b')} \ge L^{(b)}) \le \lfloor \gamma B \rfloor \right\} \end{equation} with the trimming rate $\gamma \in [0,1[$, we compute the trimmed aggregated selection frequencies \begin{equation} \label{trimself} \hat \pi_j^{\gamma}=\frac{1}{B-\lfloor \gamma B \rfloor} \sum_{b \in \{1,...,B\} \setminus I^{trim}(\gamma)} I(j \in \hat S^{(b)}) \end{equation} for all $j=1,...,p$, and let $\hat \pi^{\gamma}:=(\hat \pi_j^{\gamma})_{j=1}^p$. The actual identification of the stable set works as usual, based on the $\hat \pi_j^{\gamma}$. 

The \textbf{Trimmed Stability Selection (TrimStabSel)} is described by the following algorithm: 

\begin{algorithm}
\label{TrimStabSel}
\textbf{Initialization:} Data $D \in \R^{n \times (p+1)}$, number $n_{sub}$ of instances per resample, resampling procedure, number $B$ of resamples, either threshold $\pi_{thr}$ or number $q$ of stable variables, model selection algorithm, loss function $L$, trimming rate $\gamma$\;
\For{$b=1,...,B$}{Draw a resample $D^{(b)}=(X^{(b)},Y^{(b)}) \in \R^{n_{sub} \times (p+1)}$ from $D$\;
Apply the model selection algorithm to $D^{(b)}$\;
Get a predictor set $\hat S^{(b)}$ and coefficients $\hat \beta^{(b)}$\;
Evaluate the in-sample loss $L^{(b)}$ from Eq. \ref{insample}}
Flag the $\lfloor \gamma B\rfloor$ resamples with the highest losses as outlying\;
Compute the trimmed aggregated selection frequencies as in Eq. \ref{trimself}\;
Compute the stable set according to $\pi_{thr}$ or $q$ and $\hat \pi^{\gamma}$\;
\caption{Trimmed Stability Selection}
\end{algorithm}

Now, we have to analyze the effect of this trimming procedure on the Stab-BDP. We abstain from detailing out each of the cases considered in Subsections \ref{thrss} and \ref{rankss} again but formulate a universal result which can be easily adapted to all the individual cases. 

\begin{thm} For the robustness gain of TrimStabSel with trimming level $\gamma$ in comparison with the untrimmed Stability Selection, let the number of contaminated models in the Gaussian brackets in the right hand side of the logical expressions in $P(\cdot)$-brackets in the probabilities computed in Thm.s \ref{thrstabselrob} and \ref{rankstabselrob}, part iii) and iv) of Thm. \ref{thrstabselrobcor} resp. part c) and d) of Cor. \ref{rankstabselrobcor} be denoted by $K$. Assuming that there are $k_{\gamma}$ broken models in the set of the $\lfloor \gamma B \rfloor$ trimmed models, $k_{\gamma} \in \{0,1,...,\lfloor \gamma B \rfloor\}$, for TrimStabSel, the number $K$ is replaced by $k_{\gamma}+(B-\lfloor \gamma B \rfloor) K/B$ resp. $2k_{\gamma}+B-\lfloor \gamma B \rfloor K/B$ in the upper interval bounds in Thm. \ref{rankstabselrob} and Cor. \ref{rankstabselrobcor}. \end{thm}

\begin{bew} Since the number of aggregated models decreases from $B$ to $B-\lfloor \gamma B \rfloor$, $K$ has to be multiplied by $(B-\lfloor \gamma B \rfloor)/B$ in order to take the increasing effect of the individual non-trimmed models into account. The only missing feature is the number $k_{\gamma}$ of contaminated trimmed models which increases the allowed number of sufficiently contaminated resamples by $k_{\gamma}$. In the mentioned special cases in Thm. \ref{rankstabselrob} and Cor. \ref{rankstabselrobcor}, one has to account for the factor $0.5$ before the Gaussian bracket.  \end{bew}  \vspace{-1.25cm} \begin{flushright} $_\Box $ \end{flushright} 

\begin{ex} Let, for the threshold-based Stability Selection, $\max_j(\hat \pi_j^+)-\pi_{thr}=0.2$ and let $B=100$. Then, in Thm. \ref{thrstabselrob}, the Stability Selection does break down if at least $\lceil 100 \cdot 0.2 \rceil=20$ models have broken down. If $\gamma=0.2$ so that 20 models are trimmed and if $k_{\gamma}=20$, there are 20 broken models allowed due to trimming and additional $\lceil 80 \cdot 0.2 \rceil-1=15$ models due to the inherent robustness of the Stability Selection, requiring 36 broken models for a breakdown. Note that for $k_{\gamma}=0$, i.e., only good models are trimmed, the number of sufficiently contaminated resamples for a breakdown decreases from initially 20 to 16 which would be the price for trimming good models away. \end{ex}

\begin{rem} Another robustness gain could be achieved by putting more mass on ``good'' instances during resampling. This idea has been proposed for FRB in \cite{salibian08} who propose to apply a robust MM-estimator first in order to detect bad leverage points which are downweighted during their updating scheme. A similar approach can be very beneficial for TrimStabSel by applying such a robust model selection algorithm first on the whole data set when assuming case-wise contamination resp. a cell-wise robust model selection algorithm like the Sparse Shooting S from \cite{bottmer} so that suspicious instances may be identified and downweighted in the resampling procedure. \end{rem}

\section{Simulation study} \label{simsec}

We now investigate the impact of our proposed outlier scheme on model selection and the performance of TrimStabSel.

We consider a variety of scenarios that differ by $n$, $p$, the fraction of outlying cells and the signal to noise ratio (SNR). In all scenarios, there are $s^0=5$ relevant variables and the corresponding components of the coefficient $\beta$ are i.i.d. $\mathcal{N}(4,1)$-distributed. The cells of the regressor matrix are i.i.d. $\mathcal{N}(5,1)$-distributed. In the regression settings, the responses are computed by $Y_i=X_i\beta+\epsilon_i$ with $\epsilon _i \sim \mathcal{N}(0,\sigma^2)$ i.i.d. for all $i$ where $\sigma^2$ is set so that a specific signal to noise ratio is valid. In the classification settings, we compute $X\beta$ and $\eta_i:=\exp(\overline{X_i \beta})/(1+\exp(\overline{X_i \beta}))$ where $\overline{X_i \beta}:=X_i\beta-mean(X\beta)$. The responses $Y_i$ are drawn according to $Y_i \sim Bin(1,\eta_i)$. As there is no opportunity to simulate a target SNR, we distinguish three situations by drawing $\beta_j \sim \mathcal{N}(\mu,1)$ so that a higher $|\mu|$ corresponds to a higher SNR. As for the outlier configuration, we consider different $\tilde m \le n$ and for each $\tilde m$, we randomly select $\tilde m$ rows of the regressor matrix so that in each of these rows, the value in the cells corresponding to the five relevant variables are replaced by zero.

We consider different model selection algorithms, namely, $L_2$-Boosting (\cite{bu03}, \cite{bu07}), LogitBoost (see \cite{bu}) and SLTS (\cite{alfons13}). The stable model is always derived rank-based with $q=5$. 

We evaluate the performance of the Stability Selection variants by computing the mean true positive rate (TPR) over $V$ repetitions where for each $v=1,...,V$, we generate an independent data set. Moreover, we compute the fraction of breakdowns, i.e., where no relevant variable has entered the stable model as well as the fraction of cases where the stable model is perfect, i.e., it consists only of the five relevant variables. These quantities are then plotted against the number $\tilde m$ of outlying instances in a single graphic, separately for each SNR.

\subsection{(Trimmed) Stability Selection with $L_2$-Boosting}

The scenario specifications are given in Table \ref{scen}. We consider the SNRs 1, 2 and 5 and for each SNR and each value for $\tilde m$ from the set given in Table \ref{scen}, we generate $V=1000$ independent data sets and apply all four Stability Selection variants specified in Table \ref{scen}. We use the function \texttt{glmboost} from the $\mathsf{R}$-package \texttt{mboost} (\cite{mboost}) with \texttt{family=Gaussian()}, 100 iterations and a learning rate of $0.1$. 

\begin{table}[H]
\begin{center}
\begin{tabular}{|p{0.65cm}|p{0.55cm}|p{0.6cm}|p{1.8cm}|p{0.8cm}||p{1.25cm}||p{0.75cm}|p{0.75cm}||p{0.75cm}|p{0.75cm}||p{0.75cm}|p{0.75cm}|} \hline 
Scen. & $p$ & $n$ & $\tilde m$ & $n_{sub}$ & StabSel & \multicolumn{6}{c|}{TrimStabSel} \\ \hline
& & & & & $B$ & $B$ & $\gamma$ &  $B$ & $\gamma$ &  $B$ & $\gamma$  \\ \hline
1 & 25 & 50 & \{0,1,...,9\} & 25 & 100 & 100 & 0.5 & 100 & 0.75 & 100 & 0.9 \\ \hline
2 & 50 & 100 &  \{0,1,...,9\} & 50 & 100 & 100 & 0.5 & 100 & 0.75 & 100 & 0.9 \\ \hline
3 & 200 & 200 & \{0,1,...,20\} & 100 & 100 & 100 & 0.5 & 1000 & 0.9 & 1000 & 0.95 \\ \hline
4 & 500 & 200 & \{0,1,...,15\} & 100 & 100 & 100 & 0.75 & 1000 & 0.9 & 1000 & 0.95 \\ \hline
\end{tabular}
\end{center}
\caption[Scenario specification]{Scenario specification for $L_2$-Boosting and LogitBoost as model selection algorithms} \label{scen}
\end{table}

\begin{figure}
\begin{center}
\includegraphics[width=4.5cm]{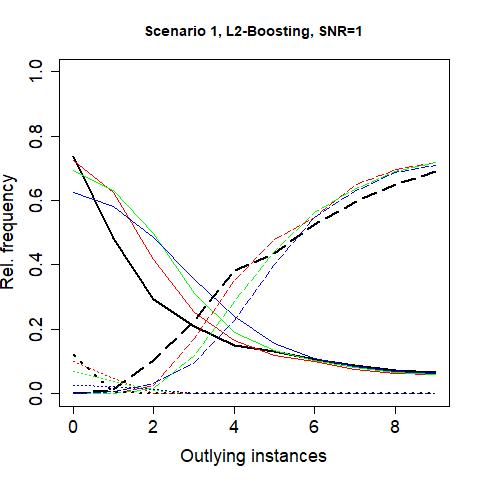} 
\includegraphics[width=4.5cm]{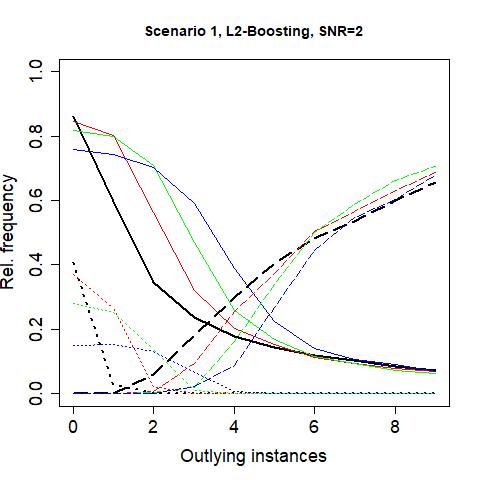} 
\includegraphics[width=4.5cm]{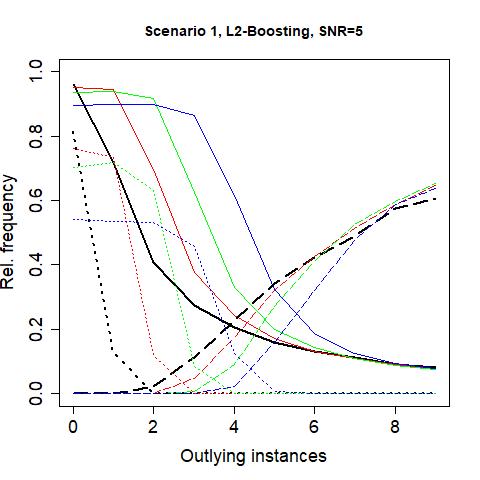} \\
\includegraphics[width=4.5cm]{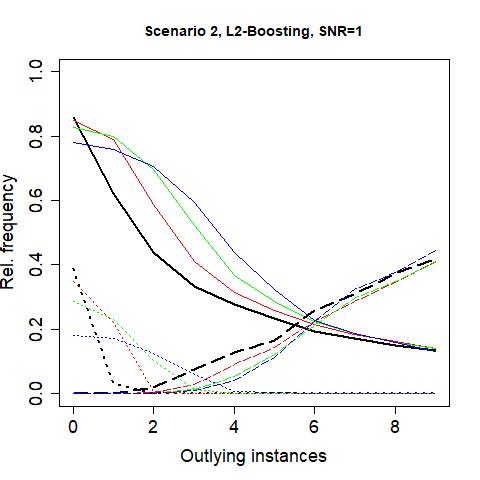} 
\includegraphics[width=4.5cm]{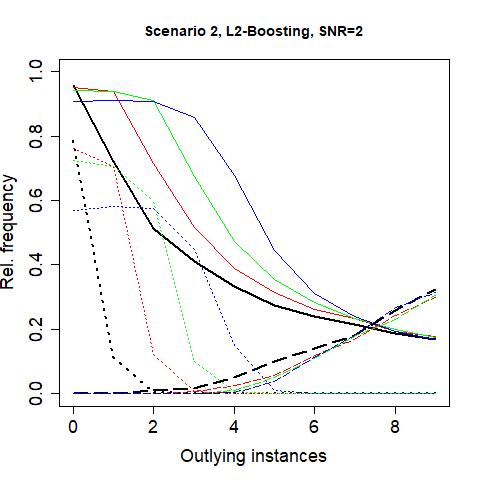}
\includegraphics[width=4.5cm]{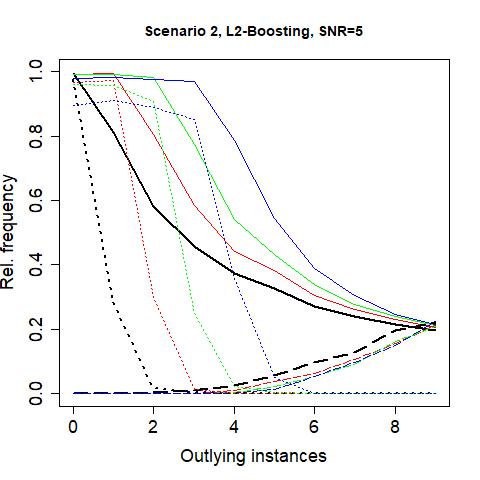}
\caption{Results for scenarios 1 and 2 with $L_2$-Boosting as model selection algorithm. Solid lines represent the TPR, dashed lines the relative frequencies of a breakdown and dotted lines the relative frequencies perfect stable models. The black lines correspond to the non-trimmed Stability Selection and the red, green and blue lines to the first, second and third configuration of TrimStabSel, as specified in Table \ref{scen}, respectively.} \label{l2results1}
\end{center}
\end{figure}

\begin{figure}
\begin{center}
\includegraphics[width=4.5cm]{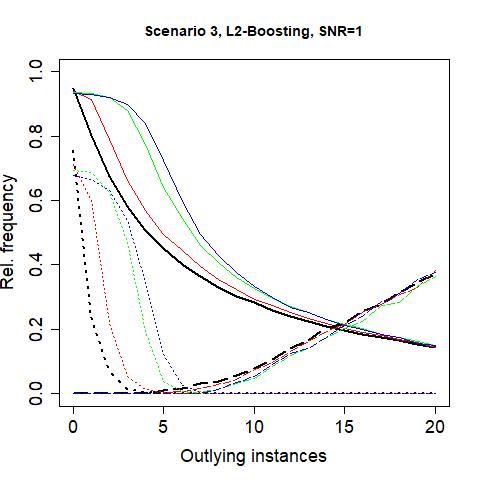} 
\includegraphics[width=4.5cm]{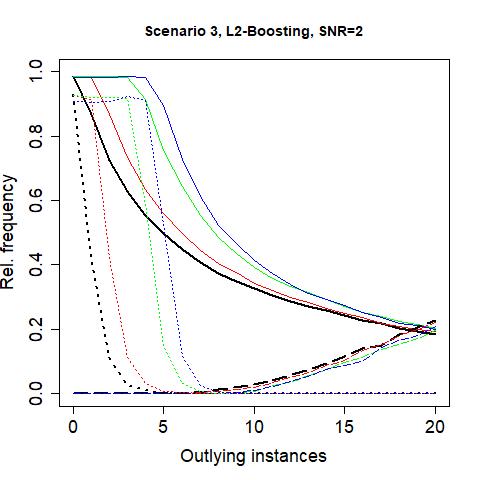} 
\includegraphics[width=4.5cm]{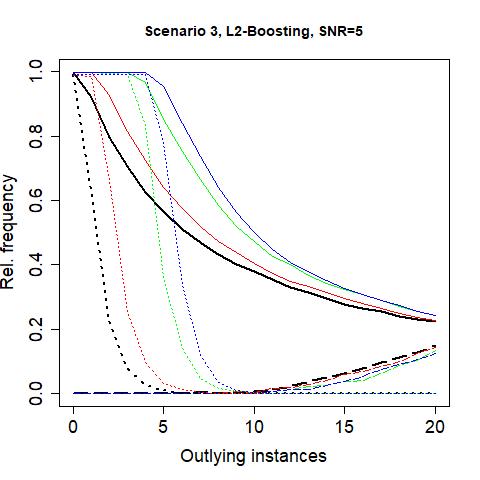} \\
\includegraphics[width=4.5cm]{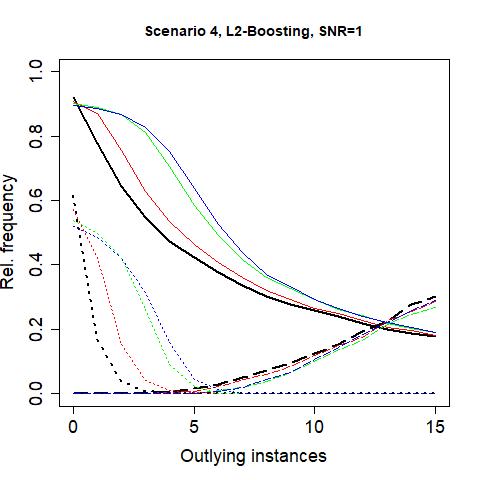} 
\includegraphics[width=4.5cm]{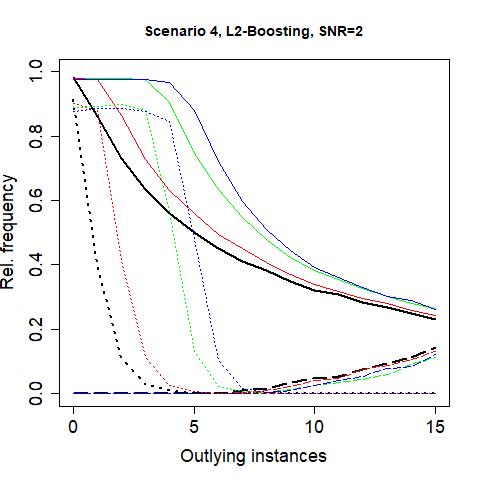}
\includegraphics[width=4.5cm]{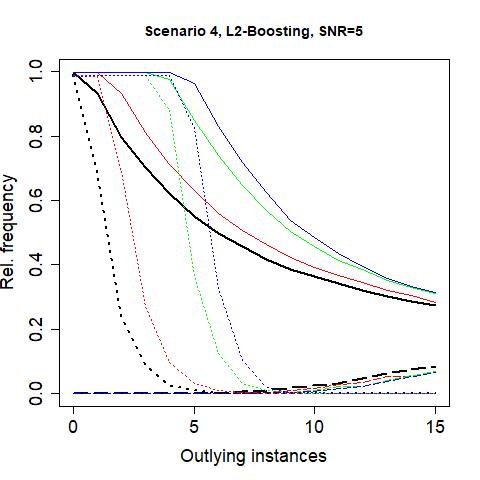}
\caption{Results for scenarios 3 and 4 with $L_2$-Boosting as model selection algorithm.} \label{l2results2}
\end{center}
\end{figure}

The results in Fig. \ref{l2results1} and Fig. \ref{l2results2} show the non-surprising facts that the TPR increases with increasing SNR and that the TPR curves and the curves representing the relative frequency of perfect models decrease with the number of outliers while the curves representing the relative frequency of a breakdown increase. For a low number of outliers, in particular for clean data, the performance of the TrimStabSel variants is usually worse than that of the non-trimmed Stability Selection due to the loss of evidence by trimming (good) models away. A characteristic aspect of all curves is that once contamination occurs, the TrimStabSel variants show better performance than the non-trimmed Stability Selection but that the robustness and performance gain decreases once too many cells are contaminated. The reason is that at some point, the expected number of contaminated subsamples becomes too high so that even TrimStabSel with the configurations in Tab. \ref{scen} is distorted. Note that the exact contamination rate where even TrimStabSel breaks down cannot be computed, but considering scenario 1, $\tilde m=6$ leads to the probability $P(Hyp(50,44,25)<25) \approx 0.989$ to draw a contaminated subsample and therefore to the probability $P(Bin(100,P(Hyp(50,44,25)<25)) \le 90) \approx 2.05 \cdot 10^{-7}$ to have at most 90 out of 100 contaminated subsamples. The inherent robustness of Stability Selection safeguards against a breakdown of the TrimStabSel variants here which would otherwise be very likely. 

At small contamination rates, for scenario 1 until around 12\% and for scenario 2 until around 7\%, TrimStabSel considerably improves model selection, especially for a high SNR. For example, the relative breakdown frequency in scenario 1 with an SNR of 5 and $\tilde m=4$ is around 10 times as high for the non-trimmed Stability Selection as for the third variant of TrimStabSel while the TPR is three times as high and even more than three times as high for $\tilde m=3$. The results in scenario 2-4 look similar as for scenario 1. For an SNR of 5, one can observe even near perfect results for at least the third TrimStabSel variant for low contamination rates up to around 3\%.

\subsection{(Trimmed) Stability Selection with LogitBoost}

We use the same specifications as in Table \ref{scen}, but as we cannot targetedly let the data have some specified SNR, we generate the relevant $\beta_j$ according to a $\mathcal{N}(\mu,1)$-distribution with $\mu \in \{1,4,8\}$ where higher means make, in expectation, the signals stronger and the SNR higher. Again, we use $V=1000$.  We again use \texttt{glmboost}, here with \texttt{family=Binomial(link='logit')} and let the other hyperparameters be as for $L_2$-Boosting.

\begin{figure}
\begin{center}
\includegraphics[width=4.5cm]{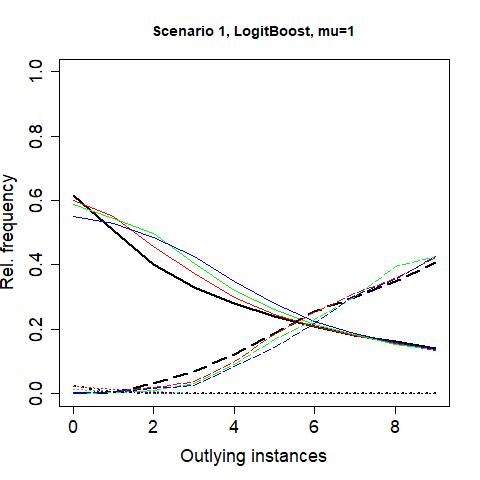} 
\includegraphics[width=4.5cm]{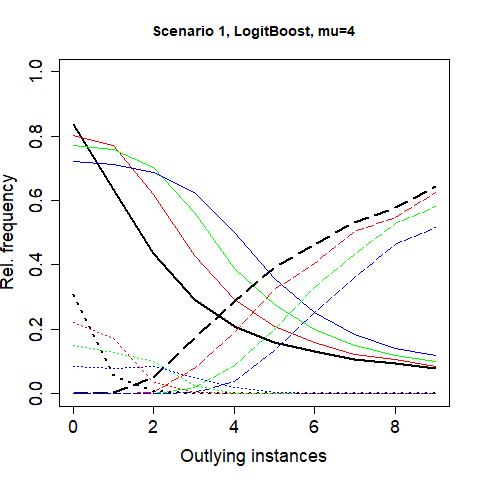} 
\includegraphics[width=4.5cm]{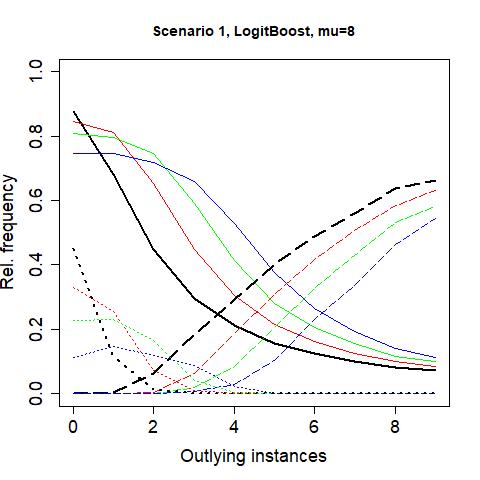} \\
\includegraphics[width=4.5cm]{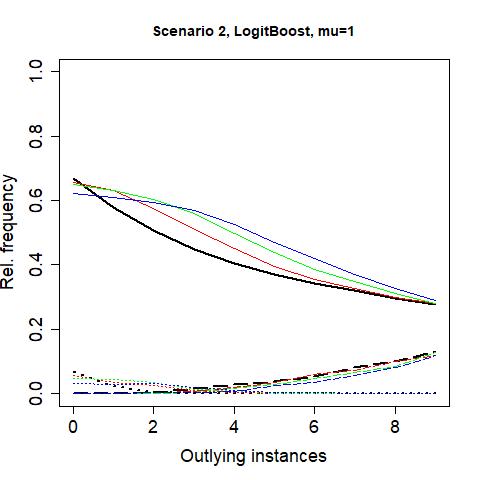} 
\includegraphics[width=4.5cm]{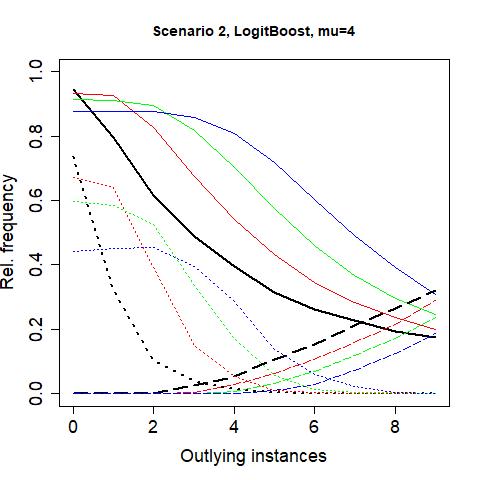}
\includegraphics[width=4.5cm]{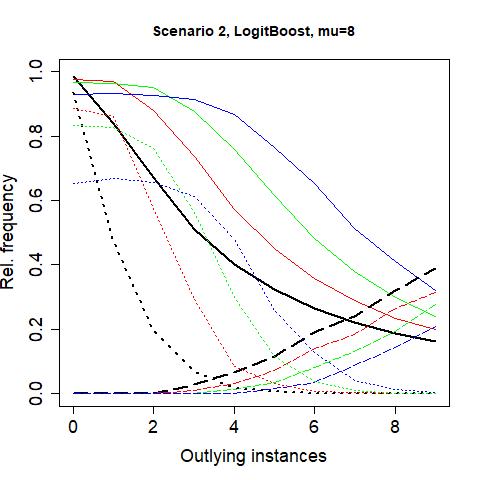}
\caption{Results for scenarios 1 and 2 with LogitBoost as model selection algorithm.} \label{glmresults1}
\end{center}
\end{figure}

\begin{figure}
\begin{center}
\includegraphics[width=4.5cm]{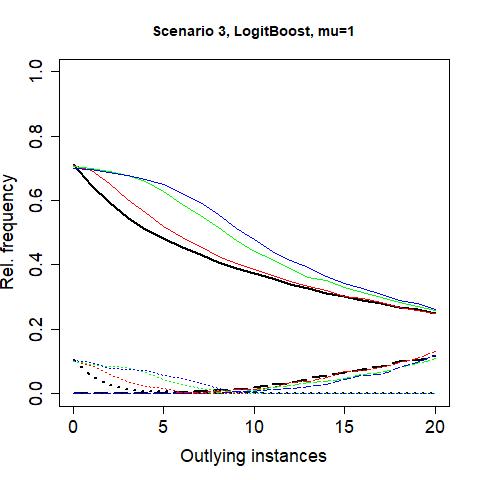} 
\includegraphics[width=4.5cm]{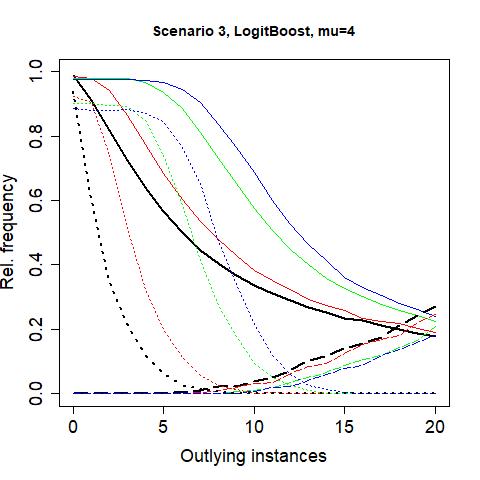} 
\includegraphics[width=4.5cm]{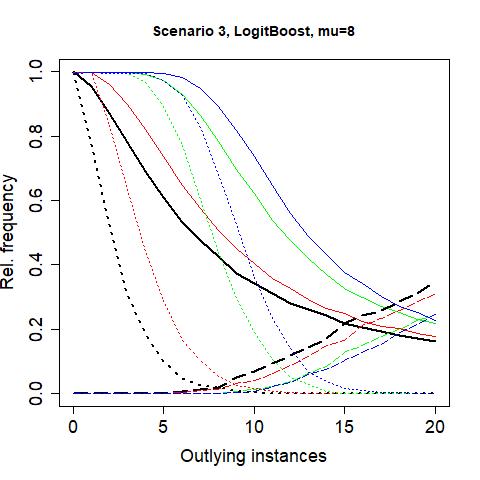} \\
\includegraphics[width=4.5cm]{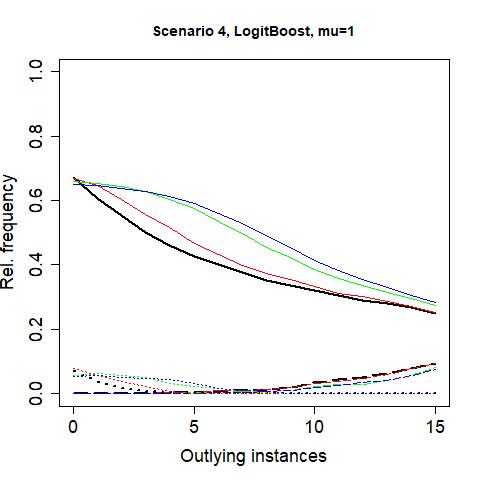} 
\includegraphics[width=4.5cm]{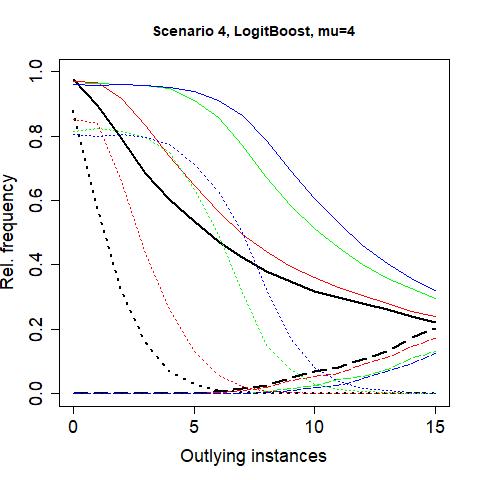}
\includegraphics[width=4.5cm]{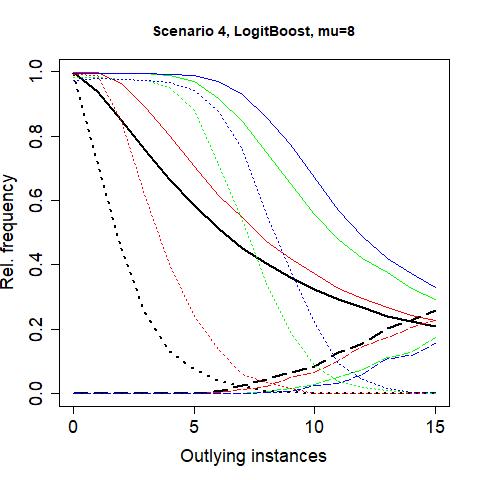}
\caption{Results for scenarios 3 and 4 with LogitBoost as model selection algorithm.} \label{glmresults2}
\end{center}
\end{figure}

The results look similarly as in the previous subsection. One can observe a more compressed shape of the TPR curve for the case $\mu=1$ while the curves corresponding to the different Stability Selection variants show a considerable margin. In our opinion, these behaviours are a direct consequence of the SNRs. The empirical mean SNRs on our data sets, computed as $Var(X\beta)/Var(Y)$ (the reciprocal value of the noise to signal ratio from \cite{elstat}), is between 30 and 40 for $\mu=1$ and more than 1200 for $\mu=8$. Although the interpretation of this SNR is not identical to the interpretation of the SNR in the regression setting, the margins between the curves for high values of $\mu$ reflect the behaviour from the previous subsection, here even stronger.

\subsection{(Trimmed) Stability Selection with SLTS}

Due to the inherent robustness of SLTS, we additionally allow for situations with high contamination radii. More precisely, we use the set $\{0,1,...,10,15,20,25\}$ for scenario 1, $\{0,1,...,10,15,20,...,45,50\}$ for scenario 2, $\{0,1,...,20,30,...,70\}$ for scenario 3 and $\{0,1,...,15,20,30,...,70\}$ for scenario 4. 

We only consider regression scenarios here, i.e., we use \texttt{family=Gaussian()} in the \texttt{sparseLTS} function from the $\mathsf{R}$-package \texttt{robustHD} (\cite{robusthd}). We use a trimming rate in SLTS of $0.25$ and for the penalty parameter, we propose the grid $\{0.05,0.55,...,4.55\}$ and let the best element be chosen data-driven by a winsorization strategy corresponding to \texttt{mode='fraction'} in \texttt{sparseLTS}. Due to the computational complexity, we set $V=100$. The remaining configurations are as in Tab. \ref{scen}.

\begin{figure}
\begin{center}
\includegraphics[width=4.5cm]{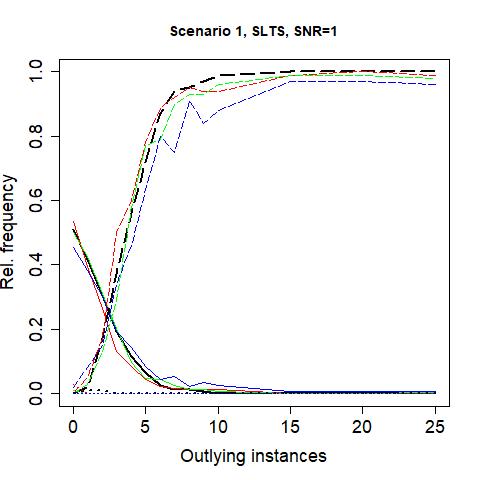} 
\includegraphics[width=4.5cm]{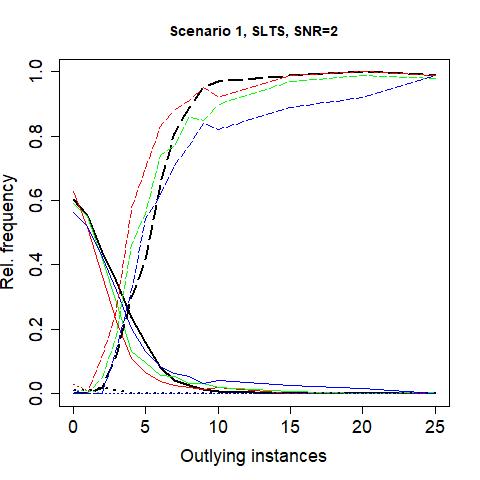} 
\includegraphics[width=4.5cm]{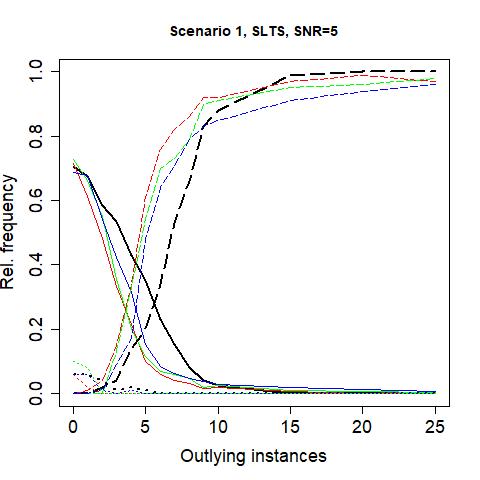} \\
\includegraphics[width=4.5cm]{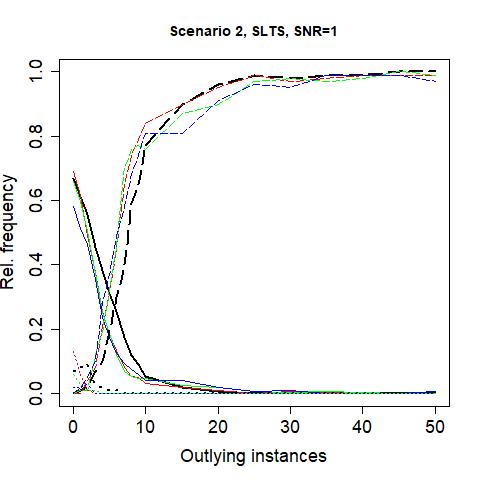} 
\includegraphics[width=4.5cm]{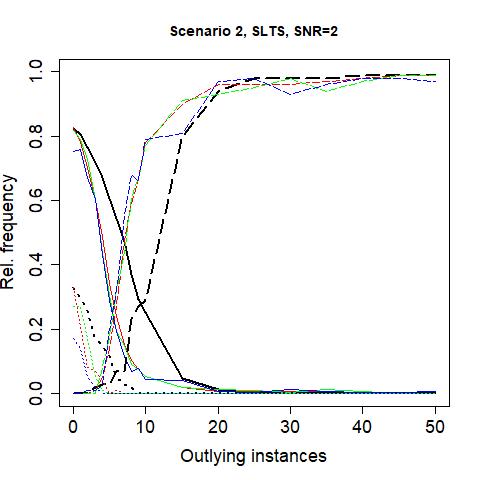}
\includegraphics[width=4.5cm]{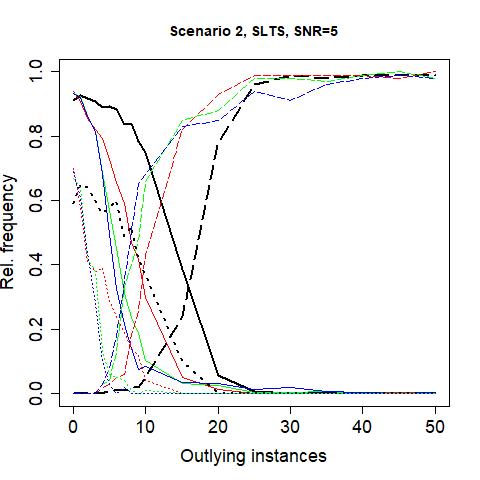}
\caption{Results for scenarios 1 and 2 with SLTS as model selection algorithm.} \label{sltsresults1}
\end{center}
\end{figure}

\begin{figure}
\begin{center}
\includegraphics[width=4.5cm]{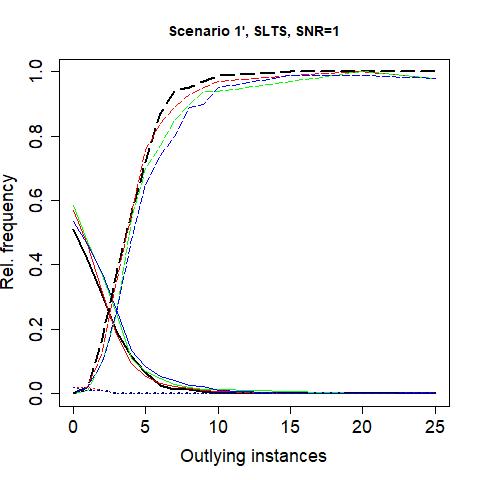} 
\includegraphics[width=4.5cm]{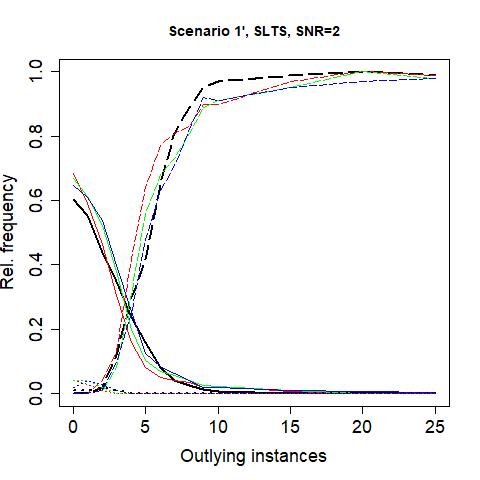} 
\includegraphics[width=4.5cm]{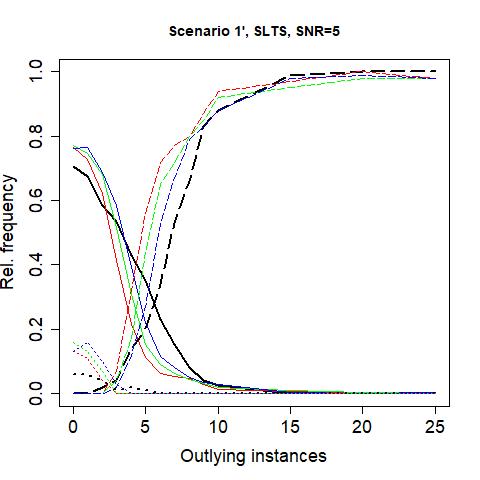} \\
\includegraphics[width=4.5cm]{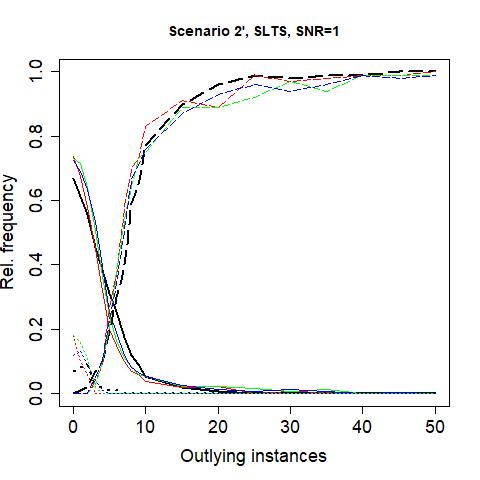} 
\includegraphics[width=4.5cm]{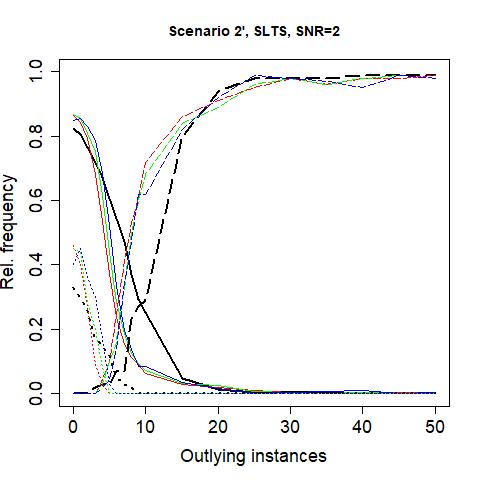} 
\includegraphics[width=4.5cm]{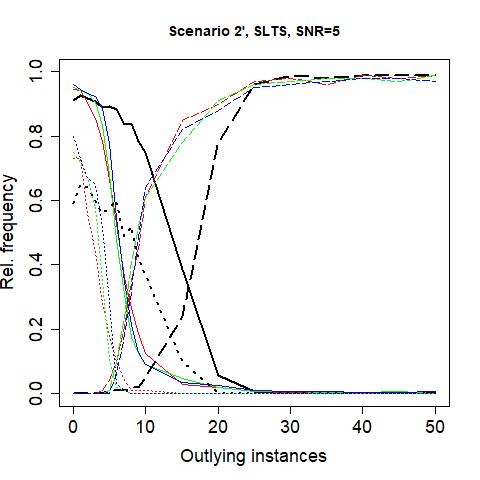} 
\caption{Results for scenarios 1' and 2' with SLTS as model selection algorithm.} \label{sltsresults2}
\end{center}
\end{figure}

\begin{figure}
\begin{center}
\includegraphics[width=4.5cm]{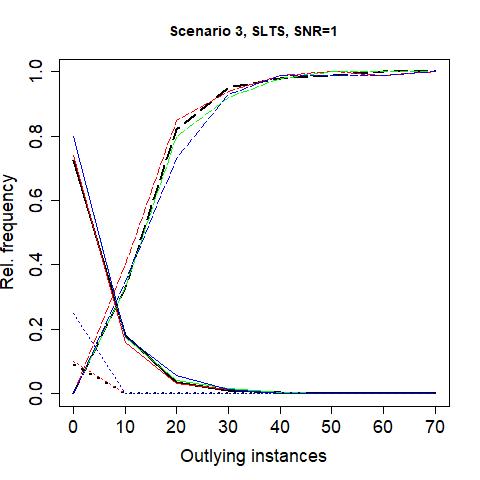} 
\includegraphics[width=4.5cm]{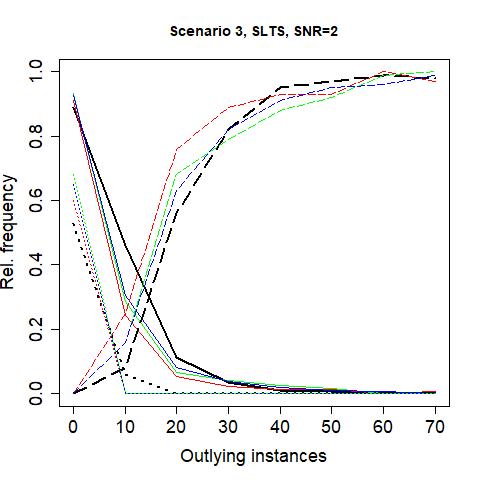} 
\includegraphics[width=4.5cm]{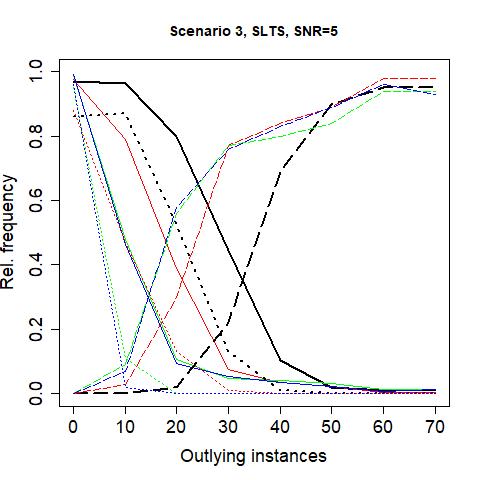} \\
\includegraphics[width=4.5cm]{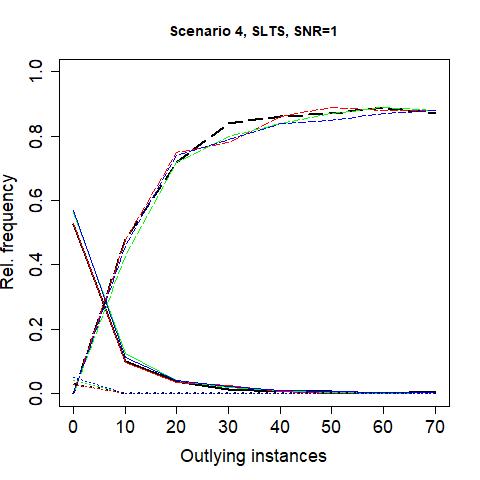} 
\includegraphics[width=4.5cm]{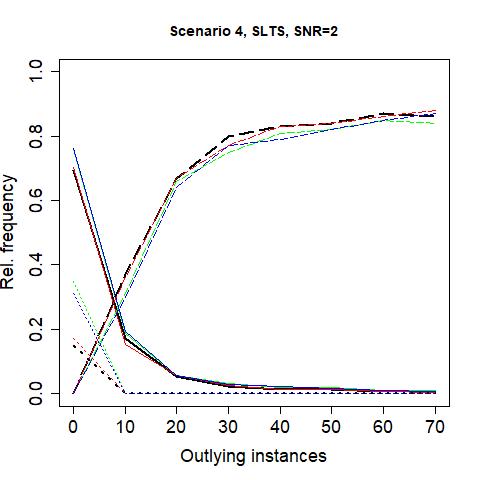} 
\includegraphics[width=4.5cm]{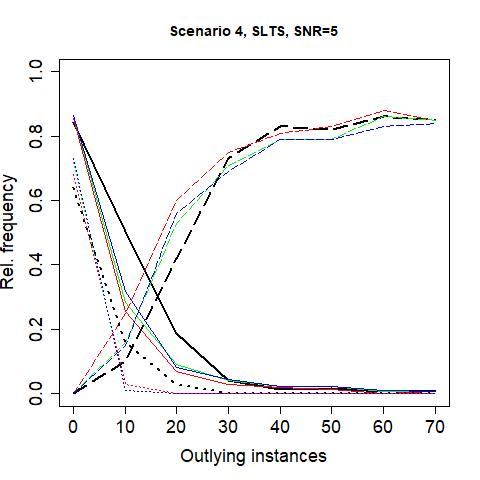} 
\caption{Results for scenarios 3 and 4 with SLTS as model selection algorithm.} \label{sltsresults3}
\end{center}
\end{figure}

The results in Tab. \ref{sltsresults1} show that the region where the non-trimmed Stability Selection leads to a better performance than the trimmed variants is considerably extended in contrast to the former experiments. This can be explained by the fact that SLTS with an internal trimming rate of 25\% itself has a BDP of around 25\% while $L_2$-Boosting and LogitBoost have a BDP of 0. Therefore, the trimmed Stability Selection pays off late, more precisely, once the contamination rate is sufficiently high to let SLTS break down on many subsamples, which is even more delayed for high SNR values due to the better performance of the underlying model selection in these cases. Also, the performance loss of TrimStabSel for low contamination rates is more considerable for low SNR values than for high SNR values. 

It should be noted that the TPR starts with lower values for $\tilde m=0$ than in the previous experiments. In particular, the relative fraction of perfect models is extremely low for scenario 1, even without contamination and with an SNR of 5, in contrast to the experiments with $L_2$-Boosting. The reason could be that, similarly as in TrimStabSel where the trimming decreases the evidence and makes the stable set therefore more fragile, trimming instances away in SLTS decreases the evidence of the fitted model and therefore its performance. The loss in efficiency of robust methods seems to carry over to model selection itself. 

One can observe that trimming eventually pays off, but due to the fraction of broken models already being very high resp. the mean TPR already being very low, the improvement itself granted by TrimStabSel may no longer be reasonable as the relative fraction of broken models is still very high, for example, in scenario 2 with an SNR of 2, 99\% of the models have broken down for the non-trimmed Stability Selection while for the third version of TrimStabSel, still 89\% of the models have broken down. 

In order to investigate the impact of the number of subsamples, we additionally propose three further TrimStabSel configurations in scenarios 1 and 2 where $B=1000$ and $\alpha \in \{0.75,0.9,0.95\}$, calling these scenarios 1' and 2'. The results are depicted in Fig. \ref{sltsresults2}. One can observe that the mean TPR is considerably higher compared to the mean TPR achieved in scenario 1 and 2, even better than for the non-trimmed Stability Selection for low contamination rates. This is a consequence of the increased evidence as we aggregate 50 up to 250 models in contrast to only 10 up to 50 models as before. Another side effect is the reduced margin between the curves corresponding to the non-trimmed Stability Selection and the TrimStabSel variants, including a smaller robustness gain in regions where the contamination rate even overcomes SLTS. The eye-catching curve of the relative breakdown rate for the third TrimStabSel variant in Fig. \ref{sltsresults1} however does no longer differ much from that of the non-trimmed Stability Selection in high contamination regions in Fig. \ref{sltsresults2}. This can be explained by solely focusing on a very few number of (non-broken) models in scenarios 1 and 2 for the third TrimStabSel variant.

\section{Conclusion}

We intended to make a step towards the unification of sparse model selection, robustness and stability in order to lift the understanding of robustness from the rows of a data matrix to the columns and investigated how contamination can affect model selection. We started with the introduction of the variable selection breakdown point and an outlier scheme which allows a very small number of contaminated cells to completely distort variable selection, making a robustification in the usual sense that provides coefficients whose norm is always bounded obsolete if no relevant variable is considered. 

We extended the notion of the resampling breakdown point, which quantifies the relative fraction of outlying instances so that the probability that a resample is contaminated too much exceeds some threshold, by the Stability Selection BDP which we computed for different scenarios where we postulate different effects of outliers onto model selection due to the absence of concrete results in literature. Our analysis reveals that a Stability Selection where the stable set is given by the best $q$ variables for a pre-defined $q$ can be expected to be more robust than the standard threshold-based Stability Selection. 

Finally, we propose a Trimmed Stability Selection which considers only the best resamples, based on the in-sample losses, when aggregating the models. A simulation study reveals the potential of this Trimmed Stability Selection to robustify model selection, although it evidently inherits the necessity to find appropriate hyperparameter configurations. The simulations also prove the alarming fragility of variable selection, even for an extremely low number of outlying cells,  if the outliers are targetedly placed onto the relevant columns. In particular, regarding the rapid performance decrease of the non-trimmed Stability Selection with $L_2$-Boosting as model selection algorithm, one has to keep in mind that even 2 resp. 5 outlying instances in scenario 1 and 2 resp. 3 and 4 suffice to let nearly no stable model be perfect, accompanied with a somewhat decreased mean TPR, so the cell-wise contamination rates range from 25/100500 in scenario 4 to 10/1275 in scenario 1. 

Although the further robustification of an SLTS-based Stability Selection in principle works which is revealed by our simulations that lead to higher mean TPRs and lower breakdown frequencies, it may not be recommended in practice. First, the trimmed Stability Selection tends to show a decreased performance in the presence of low contamination rates since SLTS can handle these configurations by its own internal trimming, similarly as a robust algorithm loses efficiency in non-contaminated settings. Second, if SLTS breaks down for large contamination rates, slightly larger contamination rates will even overcome a trimmed Stability Selection unless the number of resamples would be infeasibly high. 

We recommend to consider a trimmed Stability Selection in situations where the contamination rate can be expected to be low. In such settings, the trimmed Stability Selection with a non-robust model selection algorithm like $L_2$-Boosting shows a significant improvement concerning mean TPR, breakdown rate and the relative number of perfect stable models in contrast to the non-trimmed Stability Selection, while being very easy to implement. This avoids the application of robust model selection algorithms which are computationally more expensive and which alone do not follow the stability paradigm which would in fact necessitate to even apply a Stability Selection with a robust model selection algorithm which would be computationally very expensive. 

We want to emphasize that in our experiments, a robust model selection algorithm seems to show inferior performance than a non-robust model selection algorithm if applied on clean data. Although it is well-known that robust algorithms are less efficient in terms of asymptotic covariance, this loss in efficiency seems to carry over to variable selection itself.  

Future research is necessary in order to study the potential of outliers for targeted variable promotion or suppression further. Although our proposed outlier schemes seem to be artificial so that they most probably would not occur by chance, one has to be aware of the attacking paradigm  emerging from the deep learning community. Similarly as popular situations where models have to be inferred before attacks can be crafted (see, e.g., \cite{papernot17}), one could intercept a data transfer, try to detect relevant variables and suppress them targetedly or try to detect certainly non-relevant variables (for example by Sure Independence Screening, see \cite{fan08b}) in order to targetedly promote them.

\bibliography{Biblio}
\bibliographystyle{abbrvnat}
\setcitestyle{authoryear,open={((},close={))}}

\end{document}